\documentclass{amsart}

\usepackage{graphicx}
\usepackage{tikz}
\usepackage{stmaryrd}
\usepackage{mdframed}
\usepackage{fullpage}

\begin{document}

\title{Adaptive finite elements for obstacle problems}
\author{Tom Gustafsson$^\dag$}

\begin{abstract}
  We summarise three applications of the obstacle
  problem to membrane contact, elastoplastic torsion
  and cavitation modelling, and show how the resulting
  models can be solved using mixed finite elements.
  It is challenging to construct fixed computational
  meshes for any inequality-constrained problem
  because the coincidence set has an unknown shape.
  Consequently, we demonstrate how $h$-adaptivity can
  be used to resolve the unknown coincidence set.
  We demonstrate some practical challenges that must
  be overcome in the application
  of the adaptive method.

  \vspace{0.5cm}\noindent \emph{This preprint corresponds to the Chapter 5 of volume 58 in AAMS, Advances in Applied Mechanics.}

  \vspace{0.5cm}\noindent $^\dag$\emph{Department of Mechanical Engineering, Aalto University}
\end{abstract}


\maketitle

\section{Introduction}\label{sec1}

Let $\Omega \subset \mathbb{R}^2$ be a computational domain.
The classical obstacle problem can be written as follows: find $u : \Omega \rightarrow \mathbb{R}$ such that
\begin{equation}
  \label{eq:obs}
-\Delta u \geq f, \quad u \geq g, \quad (\Delta u + f)(u - g) = 0, \quad u|_{\partial \Omega} = 0,
\end{equation}
where $f$ and $g$ are given.
The third condition implies that in some unknown coincidence set $\Omega_C \subset \Omega$ it holds $u=g$ while
in $\Omega \setminus \Omega_C$ it holds $-\Delta u = f$.
Because $\Omega_C$ is unknown a priori,
it is challenging to construct computational meshes
without relying to mesh adaptivity.
In particular, more elements
are usually required near the interface
between the coincidence set and its complement,
and the location and the shape of the interface
are solution-dependent.

In this chapter, we review the following classical applications
of the obstacle problem:
membrane contact problem~\cite{gustafsson2017},
elastoplastic torsion problem~\cite{chouly2023},
and cavitation modelling of hydrodynamic bearings~\cite{gustafsson2018}.
We formulate the corresponding Lagrange multiplier formulation and
demonstrate numerically how a conforming $h$-adaptive and bubble-enriched
$P_2-P_0$ mixed finite
element method analysed in~\cite{gustafsson2017} performs in
each of the three applications.
Other examples of obstacle-type problems are
found, e.g., in flows through porous media~\cite{chipot2012},
and in financial mathematics~\cite{evans2012}.

Closely related but more general problems
where similar conclusions hold and, hence, analogous
conforming finite element techniques
work equally well are the frictionless contact problems
of plates~\cite{gustafsson2019}
and elastic bodies~\cite{gustafsson2020}.
However, we note that while conforming methods
have been introduced also for frictional contact problems~\cite{gustafsson2022tresca}
and some problems in plasticity~\cite{gustafsson2022bingham},
more general problems with constraints for the length of a vectorial variable
seem to be better approximated using nonconforming methods;
see also the recent monograph on the
finite element approximation of contact and friction
problems~\cite{chouly2023book}.

The method described in this
chapter is just one among many approaches for solving
obstacle problems, and has been chosen
based on its low polynomial order, simplicity
and well understood theoretical properties
such as stability,
a priori and a posteriori error estimates; cf., e.g.,~\cite{haslinger1996, veeser2001, braess2005, weiss2010, schroder2011, gustafsson2017}.
Other approaches for solving the obstacle problem
include (but are not limited to) the penalty method~\cite{gustafsson2017penalty},
Nitsche's method~\cite{chouly2023},
augmented Lagragian method~\cite{burman2017},
and the stabilised method~\cite{gustafsson2017}.

In the following. we occasionally refer to the constrained minimisation form of the obstacle problem \eqref{eq:obs}, i.e.
\begin{equation}
  \label{eq:minobs}
  \inf_{\substack{u \in H^1_0(\Omega),\\ u\geq g}} \frac12 (\nabla u, \nabla u) - (f, u)
\end{equation}
where $(.,.)$ denotes the $L^2$ inner product and $H^1_0(\Omega)$
refers to the set of functions with finite energy,
\[
(\nabla u, \nabla u) < \infty,
\]
and vanishing boundary values, $u |_{\partial \Omega} = 0$.
For more details on the related function spaces
we refer to~\cite{adams2003}.

\section{Applications}

We now review three applications where a problem of the form \eqref{eq:obs}
and \eqref{eq:minobs}
will arise.

\subsection{Deformation and contact of a membrane}

Consider a flat and thin membrane,
e.g.~$D = \Omega \times [0, t]$
where $\Omega \subset \mathbb{R}^2$
describes the shape of the membrane
on $xy$-plane and the
thickness $t>0$ is very small.
We wish to describe the out-of-plane displacement
$u(x, y) : \Omega \rightarrow \mathbb{R}$
of the membrane
under a body loading $\vec{f} = (0,0,f)$.
The governing equation is given by
\[
-\kappa \Delta u = f~\text{in $\Omega$},
\]
where $\kappa$ is a uniform (pre)tension
applied on the boundary of the membrane.

This becomes an obstacle problem if the membrane is constrained to lie above a rigid
body, $u \geq g$, whose boundary is at the distance $g$ below the membrane.
The corresponding constrained minimisation problem reads
\[
\inf_{\substack{u \in H^1_0(\Omega),\\ u \geq g}} \frac12 (\kappa \nabla u, \nabla u) - (f, u)
\]
where the membrane is fixed on the boundary.
The coincidence set $\Omega_C$ is interpreted
as the zone where the rigid body and the membrane are
in contact.
Furthermore,
the strong formulation reads: find $u : \Omega \rightarrow \mathbb{R}$ such that
\begin{equation}
  \label{eq:memobs}
-\nabla \cdot \kappa \nabla u \geq f, \quad u \geq g, \quad (-\nabla \cdot \kappa \nabla u + f)(u - g) = 0, \quad u|_{\partial \Omega} = 0,
\end{equation}

\subsection{Elastoplastic torsion problem}

\begin{figure}
  \includegraphics[width=0.7\textwidth]{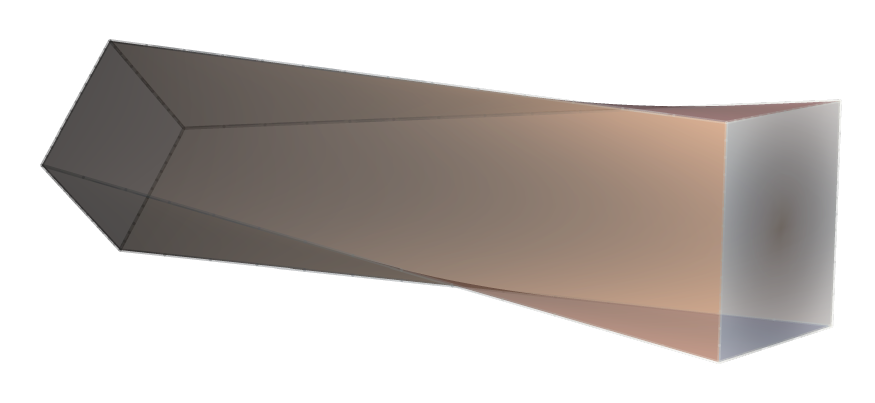}
  \caption{A visualisation of the kinematics of the torsion problem, $\vec{u}(x,y,z) = (-\theta yz, \theta xz, 0)$, when applied
    to a square shaft.  The coordinate $z$ is along the length of the shaft and $\theta$ refers to the magnitude of the rotation.
    The colouring is based on the magnitude of the displacement field
    with darker hue corresponding to smaller displacement magnitude.}
    \label{fig:torsion}
\end{figure}

Consider a uniform shaft $D = \Omega \times [0, L]$ where
$\Omega \subset \mathbb{R}^2$ is the shape of the cross section
and $L$ is the length of the shaft.
The kinematic assumption $\vec{u}(x,y,z) = (-\theta yz, \theta xz, 0)$, $0 < z < L$, $\theta > 0$,
describes a (small) rotation on $xy$-plane
whose magnitude varies linearly along the shaft while the twist per unit length $\theta$ is a given constant
and has units radians/meter, see Figure~\ref{fig:torsion}.
A similar derivation of the torsion problem
without plasticity can
be found in Feng--Shi~\cite{feng1996}.

Evaluating the infinitesimal strain and the stress tensor reveals the nonzero
stress components
\[
\sigma_{zx} = -G\theta y, \quad \sigma_{zy} = G\theta x.
\]
Differentiating with respect to $y$ and $x$, respectively,
and summing the two formulae leads to
\[
-\partial_y \sigma_{zx} + \partial_x \sigma_{zy} = 2G\theta.
\]
By introducing the stress function $\phi$ so that
\[
\partial_y \phi = \sigma_{zx} \quad \text{and} \quad -\partial_x \phi = \sigma_{zy},
\]
the above equation becomes
\[
-\Delta \phi = 2G\theta~\text{in $\Omega$}.
\]

The shaft is assumed to be traction free on the boundary of the cross section.
This is written in terms of stress as
\[
\sigma_{zx} n_x + \sigma_{zy} n_y = 0
\]
and, consequently, in terms of the stress function as
\[
\partial_y \phi n_x - \partial_x \phi n_y = 0,
\]
i.e.~the tangential derivative of $\phi$ is zero along
the boundary.  This implies that $\phi$ is constant
over the boundary but since we are only interested
in the derivatives of $\phi$ and not its value,
we can simply set $\phi|_{\partial \Omega} = 0$.

In reality there is an upper bound for the stress.
For example, the von Mises yield criterion becomes
\[
\sqrt{3(\sigma_{zx}^2 + \sigma_{zy}^2)} \leq \tau,
\]
or, in terms of the stress function,
\begin{equation}
  \label{eq:vonmises}
| \nabla \phi | \leq \frac{\tau}{\sqrt{3}},
\end{equation}
where $\tau>0$ is a given yield stress.

While \eqref{eq:vonmises} is an upper bound for
the gradient of the stress function,
Breziz--Sibony~\cite{brezis1971} show the nontrivial result that the minimisation problem
\[
\inf_{\substack{\phi \in H^1_0(\Omega),\\ |\nabla \phi| \leq \frac{\tau}{\sqrt{3}}}} \frac12 (\nabla \phi, \nabla \phi) - (2G\theta, \phi)
\]
can be reformulated as
\[
\inf_{\substack{\phi \in H^1_0(\Omega),\\ \phi \leq \frac{\tau\delta}{\sqrt{3}}}} \frac12 (\nabla \phi, \nabla \phi) - (2G\theta, \phi)
\]
where $\delta : \Omega \rightarrow \mathbb{R}$ is the distance to the boundary, i.e.
\[
\delta(\vec{x}) = \inf_{\vec{y} \in \partial \Omega} | \vec{x} - \vec{y}|.
\]
The coincidence set $\Omega_C$ is interpreted as the
zone where plastic deformation will occur.
The bound for $\phi$ is from above, but if we define $\psi = -\phi$,
we may instead consider
the following strong formulation: find $\psi : \Omega \rightarrow \mathbb{R}$ such that
\begin{equation}
  \label{eq:torobs}
-\Delta \psi \geq -2G\theta, \quad \psi \geq -\frac{\tau \delta}{\sqrt{3}}, \quad (\Delta \psi - 2G\theta)\left(\psi + \frac{\tau \delta}{\sqrt{3}}\right) = 0, \quad \psi|_{\partial \Omega} = 0.
\end{equation}

\subsection{Hydrodynamic lubrication with cavitation}

Hydrodynamic bearings are often modelled using the Reynolds equation~\cite{reynolds1886}.
The Reynolds equation for the pressure $p$ within a thin film
of lubricant,
situated between two approximately parallel surfaces,
one moving at the velocity $V$ and
at the distance $d$ while the other is fixed, reads
\begin{equation}
\label{eq:reynolds}
\nabla \cdot \frac{d^3}{\mu} \nabla p = 6 V \frac{\partial d}{\partial x},
\end{equation}
where $\mu$ is the dynamic viscosity of the lubricant.
The derivation is presented in
several references; cf.~\cite{reynolds1886, szeri2010, gustafsson2015}.
The assumptions of the model are:
\begin{enumerate}
\item The lubricant is incompressible and Newtonian.
\item The fluid film thickness $d$ is small with respect to the characteristic length
  along other dimensions.
\item The curvature of the surfaces is negligible.
\item The flow is slow enough so that the inertia can be ignored.
\end{enumerate}

\begin{figure}[h!]
  \centering
  \tikzset{
    persp/.style = {scale=1.1,x={(-0.8cm,-0.4cm)},y={(0.8cm,-0.4cm)},z={(0cm,1cm)}}
    }
    \begin{tikzpicture}[persp]
      \def\width{3}
      \def\multip{1.5}
      \def\diff{0.25}
      \fill[red!50] (0,1,0)
      \foreach \t in {5,10,...,360}
      {--(0,{cos(\t)},{sin(\t)})}--cycle;
      \fill[black!5] (0,\multip,\diff)
      \foreach \t in {5,10,...,360}
      {--(0,{\multip*cos(\t)},{\multip*sin(\t)+\diff})}--cycle;
      \fill[blue!15,draw=blue,dashed] (0,0,{\multip+\diff-2.25-0.5}) -- (0,0,1-2.25) -- (\width,0,1-2.25) -- (\width,0,{\multip+\diff-2.25-0.5}) -- cycle;

      \draw[red!50] (0,1,0)
      \foreach \t in {5,10,...,360}
      {--(0,{cos(\t)},{sin(\t)})}--cycle;
      \draw[black!50] (0,\multip,\diff)
      \foreach \t in {5,10,...,360}
      {--(0,{\multip*cos(\t)},{\multip*sin(\t)+\diff})}--cycle;

      \draw[red] (\width,1,0)
      \foreach \t in {5,10,...,360}
      {--(\width,{cos(\t)},{sin(\t)})}--cycle;
      \draw (\width,\multip,\diff)
      \foreach \t in {5,10,...,360}
      {--(\width,{\multip*cos(\t)},{\multip*sin(\t)+\diff})}--cycle;

      \def\axemulp{0.9}


      \def\fdiff{10}
      \foreach \t in {10,20,...,360} 
          \fill[red,opacity=0.2] (0,{cos(\t)},{sin(\t)}) -- (\width,{cos(\t)},{sin(\t)}) -- (\width,{cos(\t+\fdiff)},{sin(\t+\fdiff)}) -- (0,{cos(\t+\fdiff)},{sin(\t+\fdiff)}) -- cycle;
      \foreach \t in {10,20,...,360} 
          \fill[black,opacity=0.1] (0,{\multip*cos(\t)},{\multip*sin(\t)+\diff}) -- (\width,{\multip*cos(\t)},{\multip*sin(\t)+\diff}) -- (\width,{\multip*cos(\t+\fdiff)},{\multip*sin(\t+\fdiff)+\diff}) -- (0,{\multip*cos(\t+\fdiff)},{\multip*sin(\t+\fdiff)+\diff}) -- cycle;

    \end{tikzpicture}
    \hspace{0.05cm}
    \begin{tikzpicture}[scale=2.5]
        \fill[red!20] (0,0)--(1.5,0)--(1.5,1.2)--(0,1.2);
        \draw[red] (0,0)--(1.5,0);
        \draw[blue,dashed] (0,0)--(0,1.2) node[above] {$x=-R \pi$};
        \draw[blue,dashed] (1.5,0)--(1.5,1.2) node[above] {$x=R \pi$};
        \draw[red] (1.5,1.2)--(0,1.2);
        \draw (0,-0.2);
        \draw (0.4, 0.65) node [below] {$\Omega$};
        \draw[-latex,black!70] (0.75,0.02) -- (1.25,0.02) node[above] {$x$};                                                                                                                                                                                
        \draw[-latex,black!70] (0.75,0.02) -- (0.75,0.51) node[above] {$y$};
    \end{tikzpicture}
    \\[0.2cm]
    \begin{tikzpicture}[scale=0.77]
        \draw (0,0) circle (3);
        \draw[red] (0,-0.5) circle (2);
        \draw[latex-latex] (0,1.5) -- (0,2.2) node[right] {$d(x)$} -- (0,3);
        \draw[blue,dashed] (0,-3) -- (0,-2.5);
        \draw[black,domain=4:10,smooth] plot(\x,{1+0.4*cos(deg(\x-0.75))});
        \draw[red,domain=4:10,smooth] plot(\x,-1);
        \draw[blue,dashed] (4,-1) -- (4,0.4);
        \draw[blue,dashed] (10,-1) -- (10,0.4);
        \draw[-latex,black!70] (7.1,-0.9) -- (9,-0.9) node[above] {$x$};
        \draw[latex-latex] (7,-1) -- (7,0.2) node[right] {$d(x)$} -- (7,1.4);
    \end{tikzpicture}
    \caption{A schematic of a hydrodynamic bearing (top left), and its
        unraveled planar representation in a computational domain $\Omega = (-R\pi, R\pi) \times (0, L)$
        (top right). The red cylinder represents a shaft which is fitted
        inside the bearing represented by the gray cylinder. The space between
        the shaft and the bearing is filled by a thin layer of lubricant
        with thickness $d(x)$.
        The lubricant pressure is periodic over the dashed parts of the boundaries.
    }
    \label{fig:journal}
\end{figure}
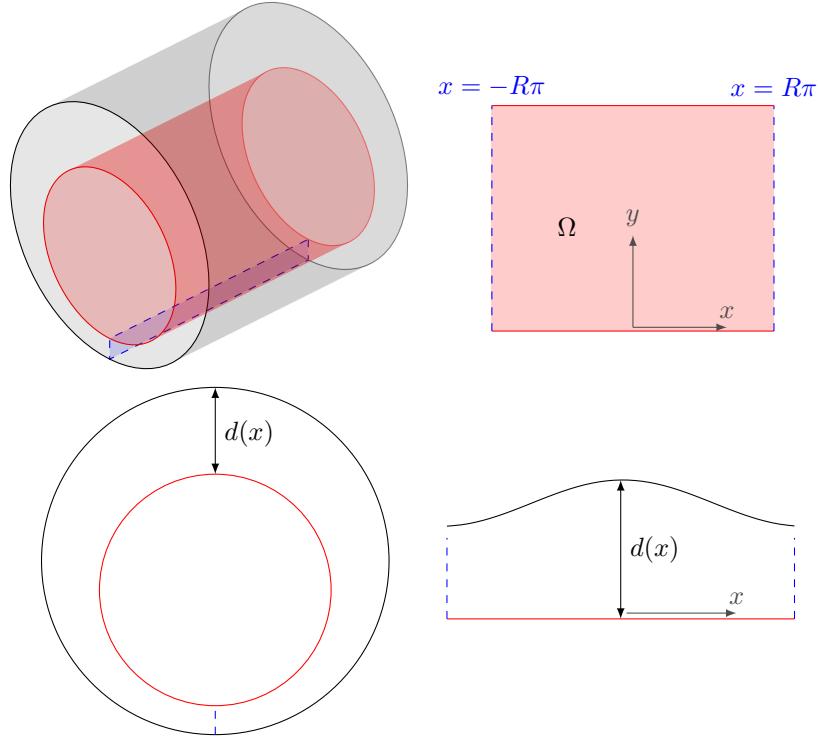

The geometry of an ideal hydrodynamic bearing, i.e.~the gap
between the housing and the shaft where the lubricant is situated, is
defined by the equation
\begin{equation}
\label{eq:gap}
d(\theta) = c_1(1 + e \cos(\theta)), \quad -\pi < \theta < \pi,
\end{equation}
where $c_1$ is the bearing clearance,
$e$ is the eccentricity, $0 \leq e < 1$.
We wish to solve \eqref{eq:reynolds}
on a rectangular domain $\Omega = (-R\pi, R\pi) \times (0, L)$
with periodic boundary conditions in $x = R\theta$,
where $R$ is the radius of the bearing
and $L$ is the width of the bearing, see Figure~\ref{fig:journal}.

Solving Reynolds equation \eqref{eq:reynolds} with
the bearing geometry \eqref{eq:gap},
boundary conditions $p(-R\pi, y) = p(R\pi, y)$, $p(x, 0) = p(x, L) = p_{\text{env}}$ and
$e \gg 0$ will result in
subatmospheric pressures on the divergent
region of the bearing.
The zone of low pressure
is subject to \emph{cavitation} which
means that any gas dissolved
within the lubricant may expand
or the lubricant may vaporise
and cause lubrication film rupture~\cite{braun2013}.
While the physical mechanisms behind cavitation
are complex,
one goal of theoretical cavitation modelling
is to incorporate film rupture into the Reynolds equation
in an
approximate manner.

The simplest cavitation model,
known as the \emph{Swift--Stieber cavitation inception condition}~\cite[p.~324]{braun2013},
states that 
\begin{equation}
\label{eq:swiftstieber}
\frac{\partial p}{\partial x} = \frac{\partial p}{\partial y} = 0, \quad p = p_{\text{cav}},
\end{equation}
where $p_{\text{cav}}$ is a constant known
as the \emph{cavitation pressure}.
While there exist several advanced cavitation
models that are physically more justified~\cite[p.~324--325]{braun2013},
this simple model has the benefit that it becomes an obstacle problem
which is well-posed and
for which exist several well understood numerical methods:
find $p : \Omega \rightarrow \mathbb{R}$ such that
\[
-\nabla \cdot \frac{d^3}{\mu} \nabla p \geq -6 V \frac{\partial d}{\partial x},\quad  p \geq p_{\text{cav}}, \quad \left( \nabla \cdot \frac{d^3}{\mu} \nabla p -6 V \frac{\partial d}{\partial x} \right)(p - p_{\text{cav}}) = 0,
\]
and
\[
p(-R\pi, y) = p(R\pi, y), \quad p(x, 0) = p(x, L) = p_{\text{env}}.
\]
The coincidence set $\Omega_C$ is the set of points where
cavitation is likely to occur.

\section{Variational formulation}

Due to the unknown coincidence set $\Omega_C$, inequality-constrained problems
are nonlinear and, hence, linearised to be solved iteratively.
Some traditional approaches to inequality-constrained problems
first discretise and then reintroduce inequality constraints
at the discrete level either by penalising or through
nodewise Lagrange multipliers.
These techniques may be suitable for a fixed mesh
and a constant obstacle but mesh adaptivity will require
imposing the inequality constraints at a continuous level.

We will apply the method of Lagrange multipliers due to
the following reasons:
(1) this will lead to an error estimate also for the
Lagrange multiplier which may have a physical interpretation such
as the contact force,
(2) the Lagrange multiplier, or an analogous quantity, is required in all finite element error
indicators we are aware of,
(3) we subjectively feel that
this approach has
the closest connection between the theory
and the implementation of the method,
(4) the implementation does not
require a reassembly of the Jacobian
if implemented as a primal-dual active
set method as demonstrated in Section~\ref{sec:linearisation}.
Another feasible approach, not discussed here,
is the use of Nitsche's method~\cite{gustafsson2017, chouly2023}
with the added benefit that the resulting linear system
can be made positive definite while the benefit of (4)
is lost.

Hence, by defining a new variable $\lambda = -\Delta u-f$,
we arrive at the alternative formulation: find $u : \Omega \rightarrow \mathbb{R}$ and $\lambda : \Omega \rightarrow \mathbb{R}$
such that
\[
-\Delta u - \lambda = f,\quad  \lambda \geq 0, \quad u \geq g, \quad \lambda(u - g) = 0, \quad u|_{\partial \Omega} = 0.
\]

From energy minimisation point-of-view, the new variable $\lambda$ can be interpreted
as a Lagrange multiplier for imposing the constraint $u \geq g$ a.e.~in $\Omega$.
In particular, the minimisation
problem
\[
\inf_{\substack{u \in V,\\ u\geq g}} \frac12 (\nabla u, \nabla u) - (f, u)
\]
can be alternatively written as the saddle point problem
\[
\inf_{u \in V} \sup_{\lambda \in \varLambda} \frac12 (\nabla u, \nabla u) - (f, u) - \langle \lambda, u - g \rangle
\]
where we have used the notation $\langle \cdot, \cdot \rangle : V' \times V \rightarrow \mathbb{R}$
for the duality pairing between the space $V$ and its topological dual space $V'$,
and the correct function spaces using standard notation for the Sobolev spaces read
\[
V = H^1_0(\Omega), \quad Q = H^{-1}(\Omega),\quad \varLambda = \{ \mu \in Q : \langle \mu, v \rangle \geq 0~\forall v \in V,~v\geq 0 \}.
\]

The finite element method is based on the variational formulation:
find $u \in V$ and $\lambda \in \varLambda$, such that
\[
(\nabla u, \nabla v) - \langle \lambda, v \rangle = (f, v) \quad \forall v \in V,
\]
and
\[
\langle \mu - \lambda, u - g \rangle \geq 0 \quad \forall \mu \in \varLambda.
\]

\subsection{Discretisation}

\newcommand{\Th}{\mathcal{T}}

Let $\Th$ denote a finite element mesh consisting of triangular elements $K \in \Th$,
and $V_\Th \subset V$, $Q_\Th \subset Q$ are the corresponding finite element spaces.
While in the continuous formulation $\lambda \in H^{-1}(\Omega)$,
which is a space of functionals,
if we use one of the standard finite element spaces for
discretising $\lambda$, say elementwise constants,
we end up with $\lambda_\Th \in L^2(\Omega)$.
Consequently, we can define the positive part of $Q_\Th$ as
\[
\varLambda_\Th = \{ \mu \in Q_\Th : \mu \geq 0 \} \subset \varLambda
\]
and the discrete variational formulation
becomes: find $u_\Th \in V_\Th$ and $\lambda_\Th \in \varLambda_\Th$ such that
\begin{equation}
  \label{eq:weak1}
(\nabla u_\Th, \nabla v) - ( \lambda_\Th, v ) = (f, v) \quad \forall v \in V_\Th,
\end{equation}
and
\begin{equation}
  \label{eq:weak2}
( \mu - \lambda_\Th, u_\Th - g ) \geq 0 \quad \forall \mu \in \varLambda_\Th.
\end{equation}

Although several stable options exist, we will now fix
the finite element spaces to bubble-enriched quadratic space for $u_\Th$,
\[
V_\Th = \{ v \in V : v|_K \in P_2(K) \oplus B_3(K) ~\forall K \in \mathcal{T} \},
\]
and piecewise-constant space for $\lambda_\Th$,
\[
Q_\Th = \{ \mu \in Q : \mu|_K \in P_0(K)~\forall K \in \mathcal{T} \},
\]
that are part of a family of bubble-enriched mixed methods analysed in Gustafsson--Stenberg--Videman~\cite{gustafsson2017}.
For this method, we have the following residual a posteriori error estimate
which is the fundamental justification behind the adaptive scheme
discussed in the upcoming sections:

\newtheorem{theorem}{Theorem}

\begin{theorem}[e.g., Gustafsson--Stenberg--Videman~\cite{gustafsson2017}]
  \label{thm:aposteriori}
  There exists $C>0$ such that
  \begin{align*}
    &\|u - u_\Th\|_{H^1(\Omega)}^2 + \|\lambda - \lambda_\Th\|_{H^{-1}(\Omega)}^2 \leq C \sum_{K \in \mathcal{T}} \bigg\{ \eta_K^2 + \eta_{\partial K}^2 + \eta_{C,K}^2 \bigg\},
  \end{align*}
  where
  \begin{align*}
    \eta_K^2 &= h_K^2 \| \Delta u_\Th + \lambda_\Th + f \|_{L^2(K)}^2 \\
    \eta_{\partial K}^2 &= \frac 12 h_K \| \llbracket \nabla u_\Th \cdot n \rrbracket\|_{L^2(\partial K \setminus \partial \Omega)}^2\\
    \eta_{C,K}^2 &=  \| (g-u_\Th)_+ \|_{H^1(K)}^2 + ((g - u_\Th)_+, \lambda_\Th)_K. \\
  \end{align*}
\end{theorem}

Above $n$ refers to the normal vector
and $\llbracket \nabla u_\Th \cdot n \rrbracket$
refers to the jump of the normal derivative over the element edge.
The positive part of $a \in H^1(\Omega)$ is denoted
by $(a)_+ = \max(a, 0)$.
Furthermore, we have used the following standard notation for the norms
in which the discretisation error is being measured:
\[
\| w \|_{L^2(S)}^2 = \int_S w^2\,\mathrm{d}x, \quad
\| w \|_{H^1(S)}^2 = \| w \|_{L^2(S)}^2 + \| \nabla w \|_{L^2(S)}^2, \quad S \subset \Omega.
\]
The Lagrange multiplier $\lambda$ belongs, in general, to
the dual space of $H^1_0(\Omega)$ which is equipped with the dual norm
\[
\| \mu \|_{H^{-1}(\Omega)} = \sup_{\substack{v \in H^1_0(\Omega),\\ \|v\|_{H^1(\Omega)} = 1}} \langle \mu, v \rangle.
\]
We note that the dual norm is not directly computable and
can be only estimated, e.g., through the mesh-dependent norm
\[
\| \mu_\Th \|_{H^{-1}(\Omega)}^2 \propto \sum_{K \in \Th} h_K^2 \| \mu_\Th \|_{0,K}^2.
\]
However, at least in the context
of the membrane contact problem \eqref{eq:memobs}, it can
be intuitively thought of as a measure similar to the total contact
force.
In particular, the total contact force would be equal to $\langle \mu, 1 \rangle$
while the dual norm maximises the force-like quantity $\langle \mu, v \rangle$ with
respect to a normalised weight function $v \in H^1_0(\Omega)$.

\subsection{Linearisation}

\label{sec:linearisation}

As explained before, inequality-constrained problems are nonlinear due to the unknown $\Omega_C$.
We will now present a Newton-type linearisation of \eqref{eq:weak2} which allows
finding $\Omega_C$ iteratively and efficiently.
The second inequality \eqref{eq:weak2} can be written using the $L^2$ projection,
$\pi_\Th : V_\Th \rightarrow Q_\Th$, as
\[
( \mu - \lambda_\Th, \pi_\Th(u_\Th - g)) \leq 0 \quad \forall \mu \in \varLambda_\Th,
\]
or, after adding and subtracting $\lambda_\Th$,
\[
( \mu - \lambda_\Th, \lambda_\Th - \pi_\Th(u_\Th - g) - \lambda_\Th ) \leq 0 \quad \forall \mu \in \varLambda_\Th.
\]
This particular form implies that $\lambda_\Th$ is equal
to the orthogonal projection of $\lambda_\Th - \pi_\Th(u_\Th - g)$
onto the set of positive finite element functions $\varLambda_\Th$.

This statement can be proven as follows: let $R=\lambda_\Th - \pi_\Th(u_\Th - g)$.
The above inequality then reads
\[
(\lambda_\Th - \mu, R - \lambda_\Th) \geq 0 \quad \forall \mu \in \varLambda_\Th.
\]
Applying this and the Cauchy--Schwarz inequality leads to
\[
\|R - \lambda_\Th \|_0^2 = (R - \lambda_\Th, R - \lambda_\Th) \leq (R - \mu, R - \lambda_\Th) \leq \|R - \mu\|_0 \|R - \lambda_\Th\|_0
\]
which holds for every $\mu \in \varLambda_\Th$.
Dividing by $\|R - \lambda_\Th\|_0$ reveals the definition of an orthogonal projection:
\[
\|R - \lambda_\Th \|_0 \leq \|R - \mu\|_0 \quad \forall \mu \in \varLambda_\Th.
\]

Because $\varLambda_\Th$ is defined as the set of nonnegative piecewise constant functions,
we can write explicitly the orthogonal projection as
\[
\lambda_\Th = (\lambda_\Th - \pi_\Th(u_\Th - g))_+.
\]
Furthermore, because $(a)_+$, $a \in \mathbb{R}$, evaluates to $a$ or 0
element-by-element, for each element either
\begin{equation}
  \label{eq:split}
\lambda_\Th = 0 \quad \text{or} \quad \pi_\Th(u_\Th - g) = 0.
\end{equation}
This suggests a solution algorithm where the split \eqref{eq:split}
is decided by the solution of the previous iteration.
Then what remains is to solve a linear problem
given by \eqref{eq:weak1} and \eqref{eq:split}
during each iteration.

As a conclusion, the \emph{primal-dual active set method} for solving the problem \eqref{eq:weak1} and \eqref{eq:weak2}
can be written as follows:
\begin{enumerate}
  \item Pick an initial guess for $u_\Th \in V_\Th$ and $\lambda_\Th \in Q_\Th$.
  \item Find the set of active elements $\mathcal{A} \subset \mathcal{T}$ such that $\lambda_\Th - \pi_\Th(u_\Th - g)|_K > 0$ for every $K \in \mathcal{A}$.
  \item Enforce $\lambda_\Th$ to zero in the complement, i.e.~$Q_{\mathcal{A}} = \{ \mu \in Q_\Th : \mu |_K = 0~\forall K \in \mathcal{T} \setminus \mathcal{A} \}$.
  \item Find $u_\Th \in V_\Th$ and $\lambda_\Th \in Q_{\mathcal{A}}$ such that
    \[
(\nabla u_\Th, \nabla v) - (\lambda_\Th, v) - (u_\Th, \mu) = (f, v) - (g, \mu) \quad \forall (v, \mu) \in V_\Th \times Q_{\mathcal{A}}.
    \]
  \item Check if the iteration has converged by calculating how much $\lambda_\Th$ changed.
    If the iteration did not converge, go back to step 2.
\end{enumerate}
The convergence of this iteration is rapid because it is,
in fact, a (semismooth) Newton's method in disguise~\cite{hintermuller2002}.
In particular, the 
number of iterations required to converge to
any given tolerance (say machine epsilon) is often
several order of magnitudes less than in
any Uzawa-type fixed point
methods we have attempted to use.
We have used the relative change in the $l^2$-norm
of the vector representing $\lambda_\Th$
as a convergence criterion,
see Listing~\ref{src:primaldual}
for our Python implementation.

\mdtheorem[style=BtypeA]{sourcecode}{Listing}

\newpage
\begin{sourcecode}[Primal-dual active set method using scikit-fem~\cite{gustafsson2020scikit}]
\label{src:primaldual}
\begin{verbatim}
# to install dependencies: pip install scikit-fem==9.0.0 matplotlib

from skfem import *
from skfem.helpers import *

mesh = MeshTri.init_sqsymmetric().refined(3)
basis = Basis(mesh, ElementTriP0() * ElementTriP2B())

def f(x):
    return 0 * x[0]

def g(x):
    return np.sin(np.pi * x[0]) * np.sin(np.pi * x[1]) - 0.5

@BilinearForm
def bilinf(lam, u, mu, v, _):
    return dot(grad(u), grad(v)) - lam * v - mu * u

@LinearForm
def linf(mu, v, w):
    return f(w.x) * v - g(w.x) * mu

lamu = basis.zeros()
(lam, lambasis), (u, ubasis) = basis.split(lamu)
A = bilinf.assemble(basis)
b = linf.assemble(basis)
G = lambasis.project(g)
res = 1e8

while res > 1e-10:
    (lam, lambasis), (u, ubasis) = basis.split(lamu)
    U = lambasis.project(ubasis.interpolate(u))
    D = np.concatenate((
        basis.get_dofs().all("u^2"),
        basis.get_dofs(elements=np.nonzero(U - G >= lam)[0]).all("u^1"),
    ))
    lamuprev = lamu.copy()
    lamu = solve(*condense(A, b, D=D))
    res = np.linalg.norm(lamu - lamuprev) / np.linalg.norm(lamu)
    print(res)

# to plot: ubasis.plot(u).show()
\end{verbatim}
\end{sourcecode}

\subsection{Implementation}

The main ideas behind adaptive
finite element schemes are given
thoroughly in Verf\"{u}rth~\cite{verfurthbook}.
The theory behind the convergence of the adaptive scheme
in case of obstacle problems
has been discussed in~\cite{braess2007}.
We will now briefly summarise the relevant parts
on the implementation of adaptive methods
for context and
to introduce notation.
Adaptive finite element methods require:
\begin{enumerate}
\item an initial mesh $\Th_0$
  \item an error indicator $E_K = E_K(u_{\Th}, \lambda_{\Th})$ which can be computed for every element $K \in \mathcal{T}$ given the discrete solution $(u_{\Th}, \lambda_{\Th}) \in V_{\Th} \times \varLambda_{\Th}$
  \item a marking strategy which decides on a subset of elements $\widetilde{\Th} \subset \Th$ to refine based on the values of the error indicator $E_K$, $K \in \Th$
  \item a mesh refinement algorithm which produces the next mesh $\Th_{i+1}$ given a subset of elements to refine on the previous mesh $\widetilde{\Th}_i \subset \Th_i$, $i=0,1,\dots$.
\end{enumerate}

\begin{figure}
  \includegraphics[width=0.3\textwidth]{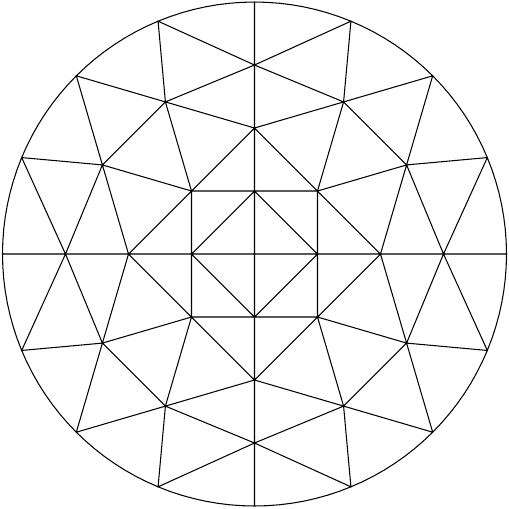}
  \includegraphics[width=0.3\textwidth]{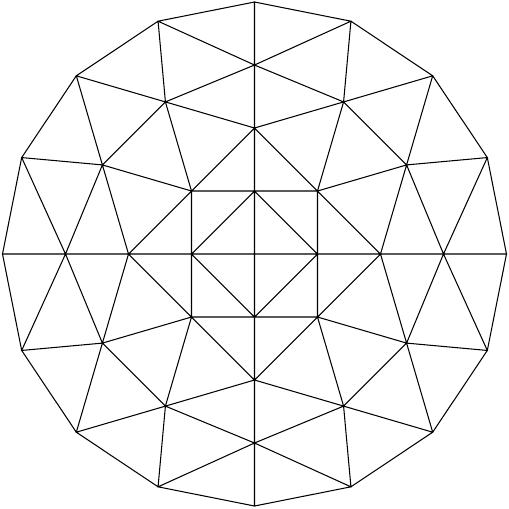}
  \caption{Adaptive schemes may overrefine if a curved domain (left) is
    approximated by a polygonal domain (right) without modifying
    the error indicator.}
  \label{fig:curvedvslinear}
\end{figure}

We wish to note that in adaptive schemes,
it is preferable to have an initial mesh $\mathcal{T}_0$ which accurately
represents the geometry of the problem because
the estimators may erroreously react to \emph{variational crimes},
e.g., replacing a circular boundary with linear segments (see Figure~\ref{fig:curvedvslinear}) can
cause excess refinement near the boundary while none might be truly required.
This is a real fundamental challenge in the practical application
of adaptive schemes: for analysis we need to assume
some explicit form for the discrete spaces
while accurately representing the geometry may require
higher order mappings and, hence, different finite element spaces
that may or may not warrant for a new analysis using different estimators.
For now we assume that the domain $\Omega$
can be specified exactly by a triangular mesh.

\begin{figure}
  \includegraphics[width=0.5\textwidth]{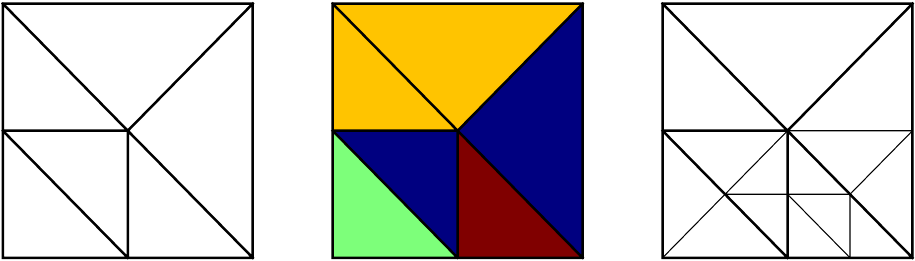}
\caption{(Red-green-blue refinement) The mesh on the left is being refined by marking one triangle for refinement.
  The red triangle is marked and split into four.
  To avoid hanging nodes and reduction in the mesh quality,
  the blue triangles are split into three
  and the green triangle is split into two.  The yellow triangles are not split.
  The refined mesh is depicted on the right.}
\label{fig:refdom}
\end{figure}

\begin{figure}
  \includegraphics[width=0.7\textwidth]{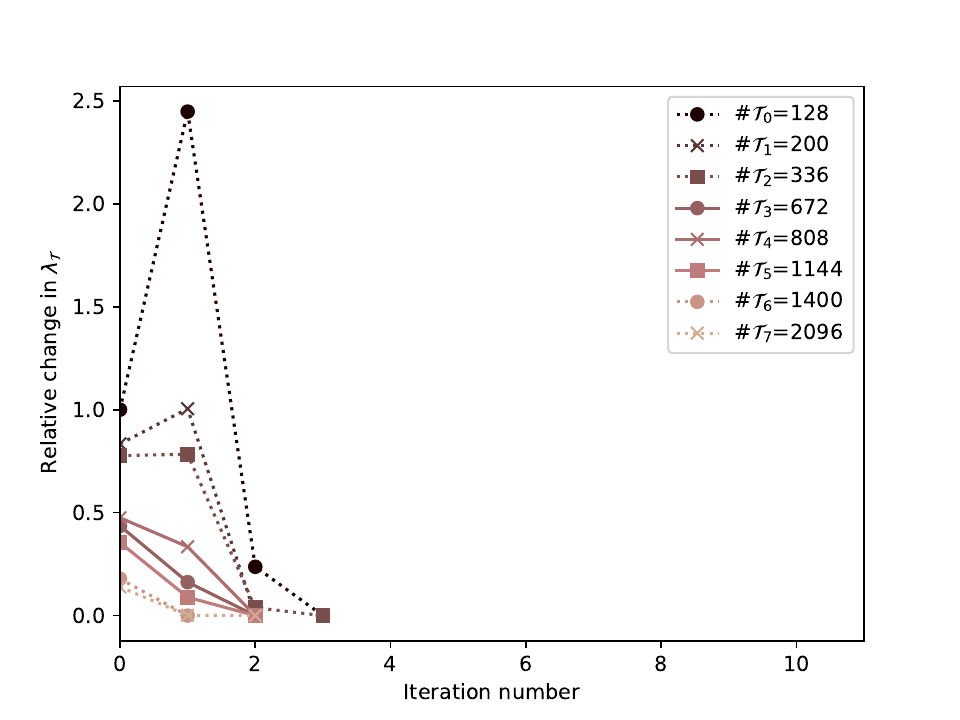}
  \includegraphics[width=0.7\textwidth]{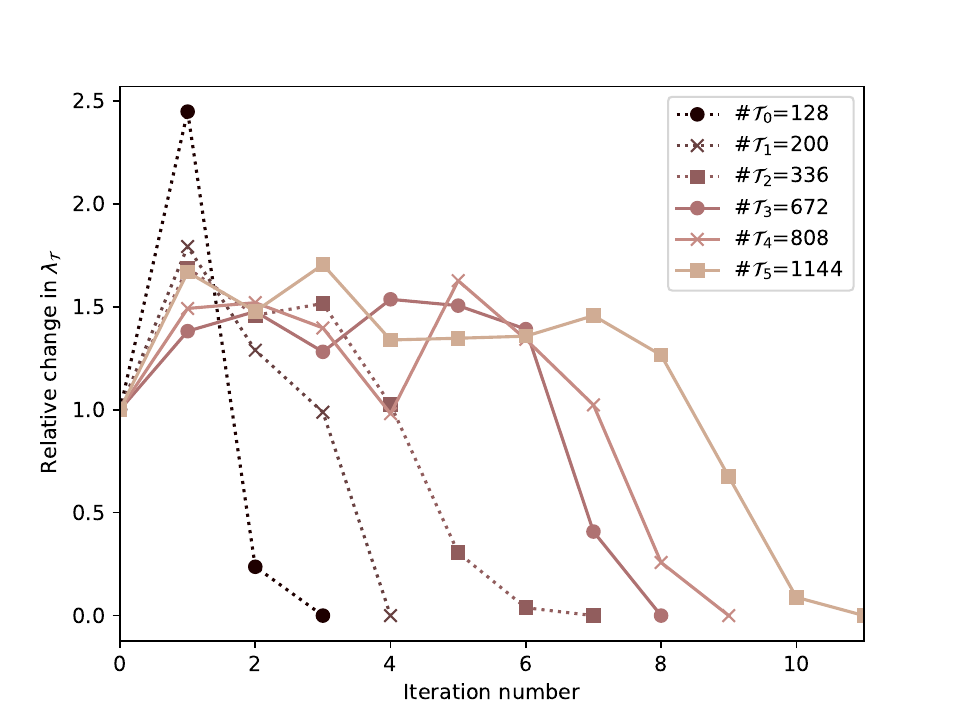}
\caption{(Primal-dual active set method) The iterations required for the primal-dual active set method
  to converge can be reduced significantly by projecting the previous solution
  to the refined mesh (top) versus starting from a zero
intial guess (bottom).}
\label{fig:primaldual}
\end{figure}

It is beneficial to have an initial mesh $\mathcal{T}_0$
which consist of shape-regular elements.  While this can
be avoided using mesh smoothing, performing refinement
without smoothing is useful because
this leads to a sequence of meshes
\[
\mathcal{T}_0 \supset \mathcal{T}_1 \supset \mathcal{T}_2 \supset \dots
\]
where every element in the previous mesh $\mathcal{T}_i$
is a union of smaller or equal elements in the subsequent mesh $\mathcal{T}_{i+1}$.
Consequently, the converged solution on $\mathcal{T}_i$
can be interpolated on $\mathcal{T}_{i+1}$ without much difficulty
and computational effort
to obtain an initial guess for the next primal-dual iteration.
Use of an initial guess from the previous mesh
will greatly reduce the number of iterations required
for the primal-dual iteration to converge; see, e.g., Figure~\ref{fig:primaldual}.

After calculating the discrete solution corresponding to the mesh $\Th_i$,
we evaluate the error
indicator on each element $K \in \Th_i$:
\begin{equation}
  E_K = \sqrt{\eta_K^2 + \eta_{\partial K}^2 + \eta_{C,K}^2}.
\end{equation}
Next we decide which elements to refine based on a marking strategy,
referred to as \emph{maximum strategy}
in Verf\"{u}rth~\cite{verfurthbook}, which marks $K \in \Th_i$
for refinement if it is larger than some given parameter $\beta \in (0, 1)$ times the
largest error indicator:
\[
\widetilde{\Th}_i = \left\{ K \in \Th_i : E_K \geq \beta \max_{K' \in \Th_i} E_{K'} \right\}.
\]
The parameter $\beta$ is used to control how many elements are refined
during each iteration of the adaptive scheme -- larger values may lead to
better final meshes and lower error (in absolute terms) but this also
drastically increases the
number of linear solves.
For simplicity, we have set $\beta = \tfrac12$
which, according to  Verf\"{u}rth~\cite{verfurthbook}, is a popular and well-established
choice.
In practice, we have observed that changing $\beta$ does not seem
to affect the asymptotic convergence rate.

Finally, we split the marked triangles using a mesh refinement algorithm,
referred to as \emph{regular red refinement}
or \emph{red-green-blue refinement},
which produces $\Th_{i+1}$ given $\widetilde{\Th}_i$.
The main idea is to split each triangle in $\widetilde{\Th}_i$ into four
while the neighboring elements are split into two or three triangles
depending on the lengths of the edges, to avoid hanging nodes
and element quality reduction; see Figure~\ref{fig:refdom}
for an example of the mesh refinement algorithm
applied to a tiny mesh.
See Listing~\ref{src:errorestim}
for our Python implementation
of the adaptive refinement.

\newpage

\begin{sourcecode}[Evaluate error estimators using scikit-fem~\cite{gustafsson2020scikit}]
\label{src:errorestim}
\begin{verbatim}
# before use: run the solver from above

gradbasis = basis.with_element(ElementDG(ElementTriP3()))
ux = gradbasis.project(ubasis.interpolate(u).grad[0])
uy = gradbasis.project(ubasis.interpolate(u).grad[1])

@Functional
def interior(w):
    return w.h ** 2 * (w["ux"].grad[0] + w["uy"].grad[1]
                       + w["lam"] + f(w.x)) ** 2

eta_int = interior.elemental(gradbasis, ux=ux, uy=uy,
                             lam=lambasis.interpolate(lam))

fbasis = [InteriorFacetBasis(mesh, ElementTriP2B(), side=i)
          for i in [0, 1]]
u12 = {"u{}".format(i + 1): fbasis[i].interpolate(u) for i in [0, 1]}
    
@Functional
def edge_jump(w):
    return w.h * dot(grad(w["u1"]) - grad(w["u2"]), w.n) ** 2

eta_edge = edge_jump.elemental(fbasis[0], **u12)    
tmp = np.zeros(mesh.facets.shape[1])
np.add.at(tmp, fbasis[0].find, eta_edge)
eta_edge = np.sum(.5 * tmp[mesh.t2f], axis=0)

def gx(x):
    return np.pi * np.cos(np.pi * x[0]) * np.sin(np.pi * x[1])

def gy(x):
    return np.pi * np.sin(np.pi * x[0]) * np.cos(np.pi * x[1])

@Functional
def contact(w):
    ind = np.maximum(g(w.x) - w["u"], 0)
    return ind * (np.maximum(g(w.x) - w["u"], 0)
                  + gx(w.x) - w["u"].grad[0]
                  + gy(w.x) - w["u"].grad[1]
                  + w["lam"])

eta_contact = contact.elemental(basis,
                                u=ubasis.interpolate(u),
                                lam=lambasis.interpolate(lam))

eta = eta_int + eta_edge + eta_contact

# to plot: mesh.plot(eta).show()
# to refine: mesh.refined(adaptive_theta(eta)).draw().show()
\end{verbatim}
\end{sourcecode}

\section{Computational demonstrations}

The following numerical results
were created using scikit-fem~\cite{gustafsson2020scikit}
which relies heavily on Python scientific computing ecosystem
\cite{harris2020array, virtanen2020scipy, hunter2007matplotlib}.
All results are obtained using minor variations of the
source code in Listings~\ref{src:primaldual} and \ref{src:errorestim}.
The mesh generation for the torsion example was done
using a variant of distmesh~\cite{persson2004simple}.

\subsection{Membrane contact problem}

\begin{figure}
  \includegraphics[width=0.32\textwidth]{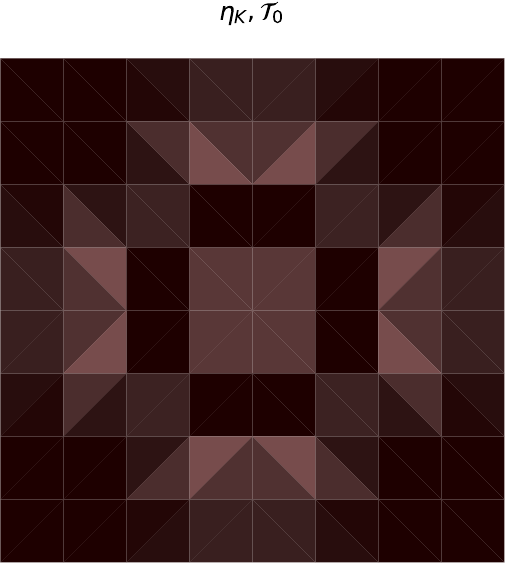}
  \includegraphics[width=0.32\textwidth]{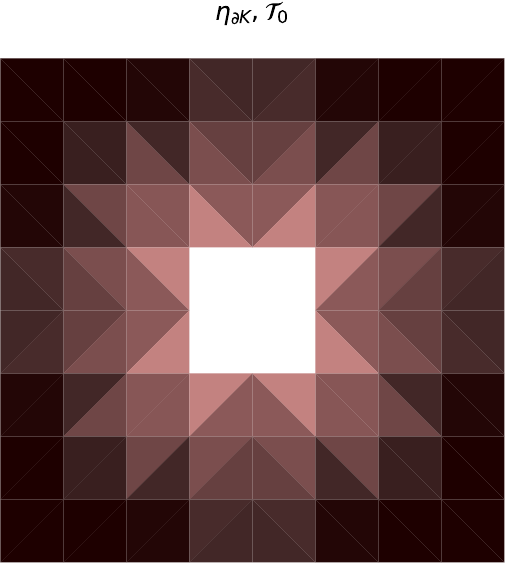}
  \includegraphics[width=0.32\textwidth]{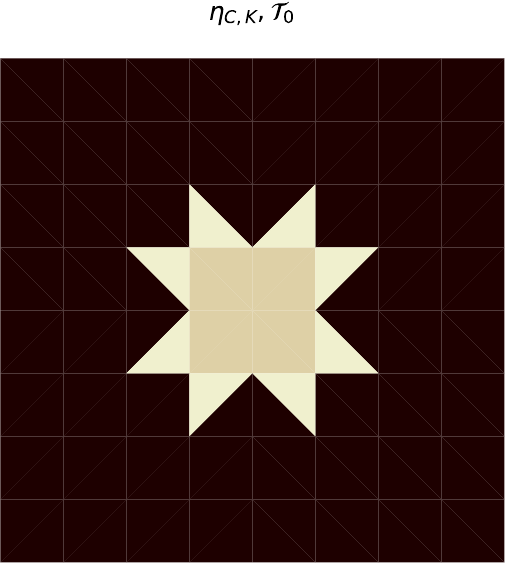}\\[0.2cm]
  \includegraphics[width=0.32\textwidth]{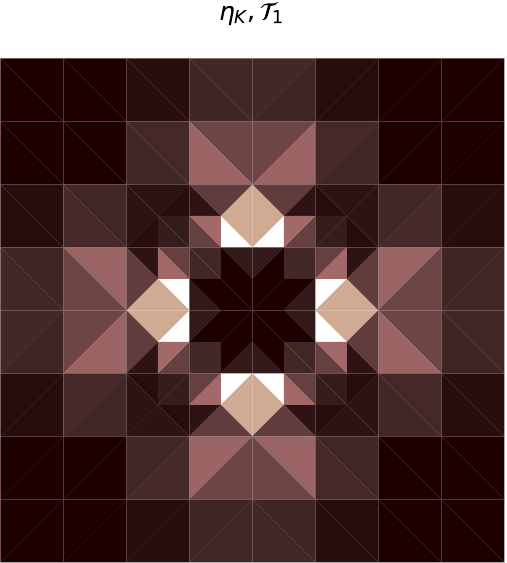}
  \includegraphics[width=0.32\textwidth]{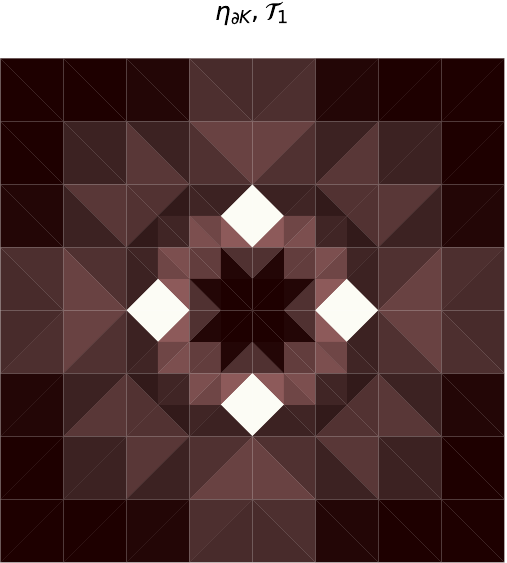}
  \includegraphics[width=0.32\textwidth]{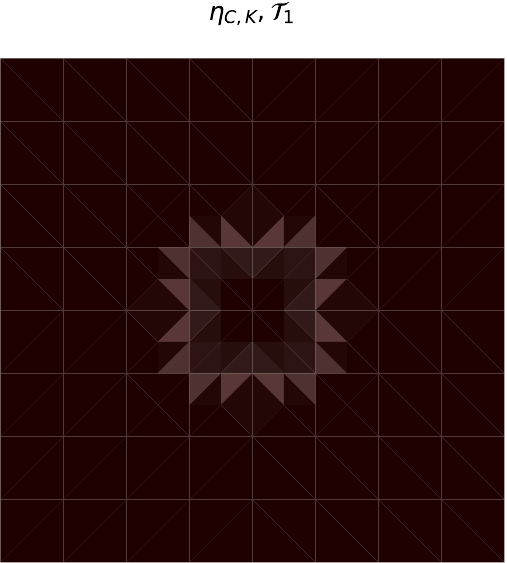}

\caption{(The sequence of error indicators for the membrane problem, part 1/3) The same colourmap is used for each indicator.
  Their exact values are of no importance as the mesh refinement
  is driven by their relative differences.
  In particular, elements within a 50 \% difference
  of the maximum value are split during each iteration.}
\label{fig:membraneetas}
\end{figure}

\begin{figure}
  \includegraphics[width=0.33\textwidth]{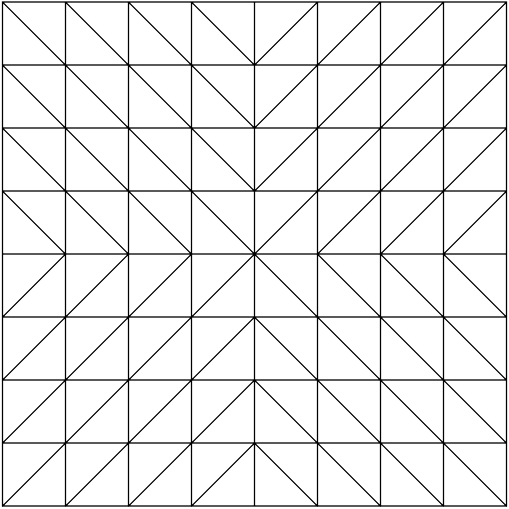}
  \includegraphics[width=0.33\textwidth]{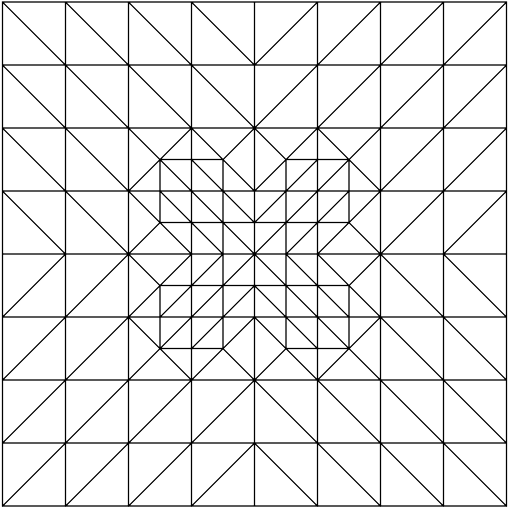}
  \includegraphics[width=0.33\textwidth]{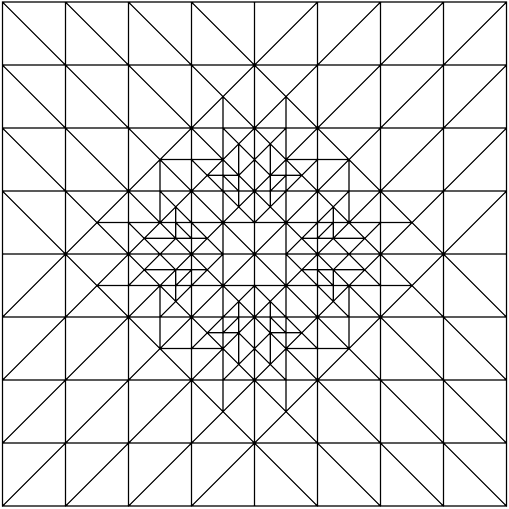}
  \includegraphics[width=0.33\textwidth]{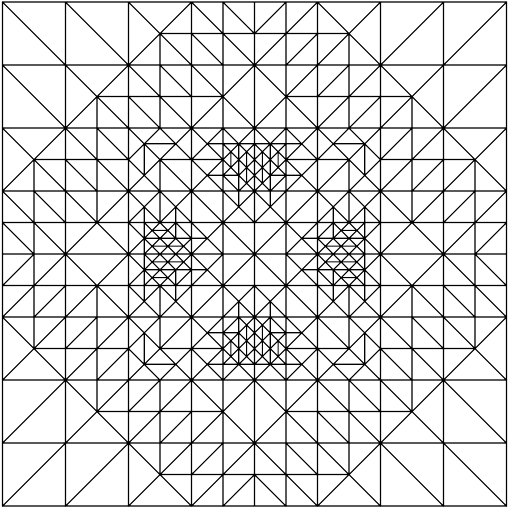}
  \includegraphics[width=0.33\textwidth]{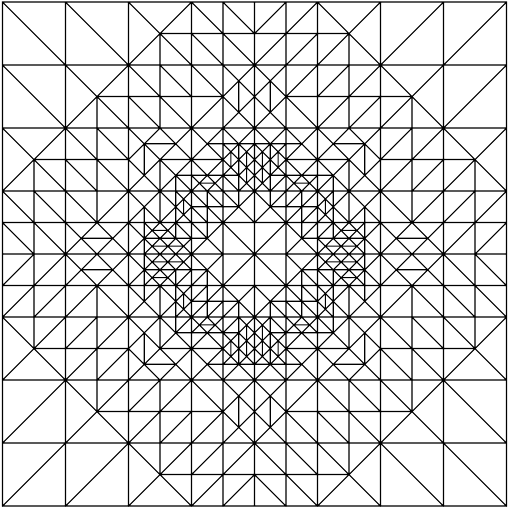}
  \includegraphics[width=0.33\textwidth]{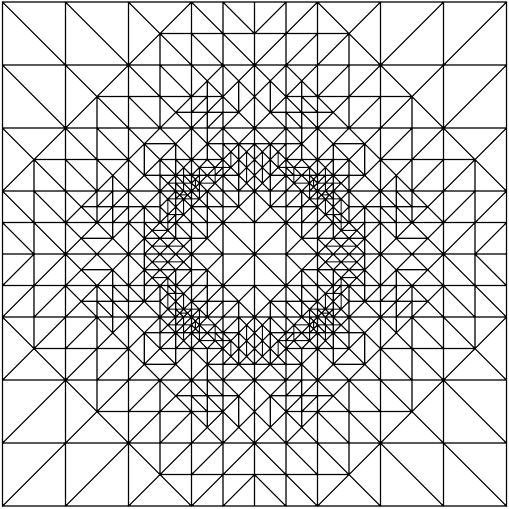}
\caption{(The sequence of adaptive meshes for the membrane problem)}
\label{fig:membranemeshes}
\end{figure}

\begin{figure}
  \includegraphics[width=0.33\textwidth]{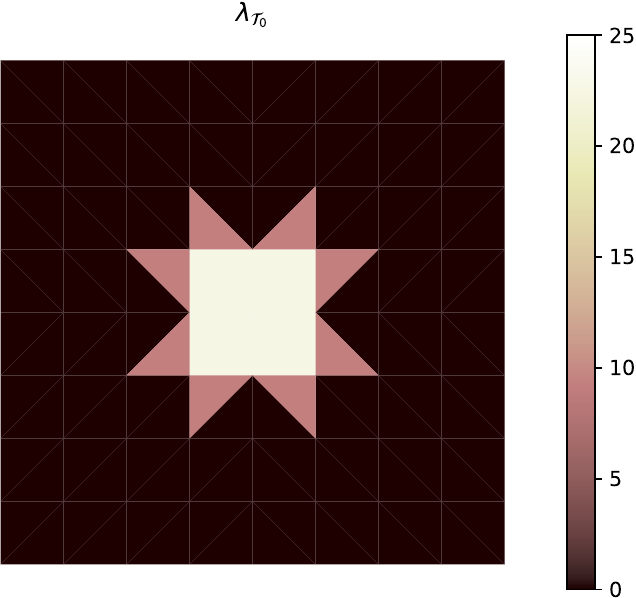}
  \includegraphics[width=0.33\textwidth]{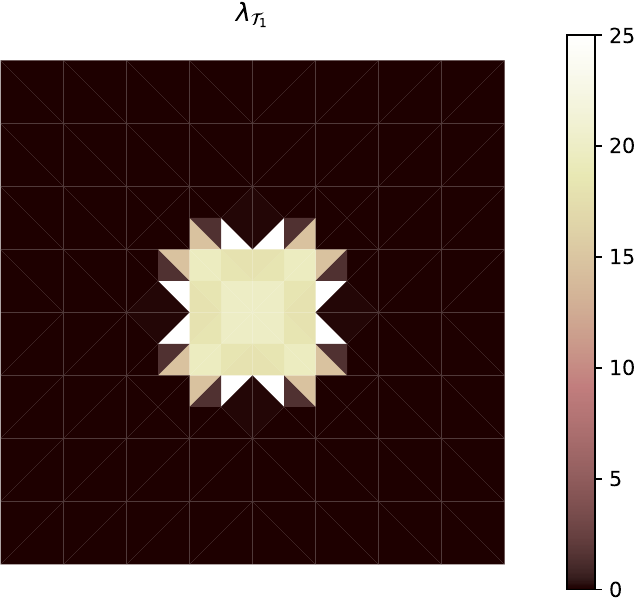}
  \includegraphics[width=0.33\textwidth]{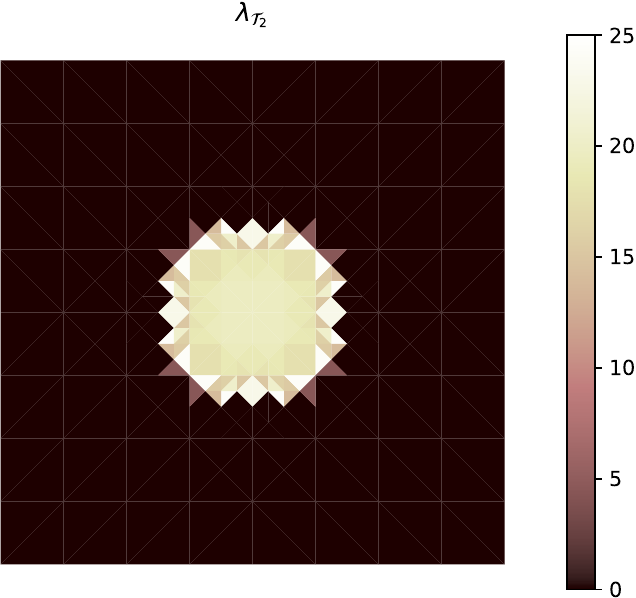}
  \includegraphics[width=0.33\textwidth]{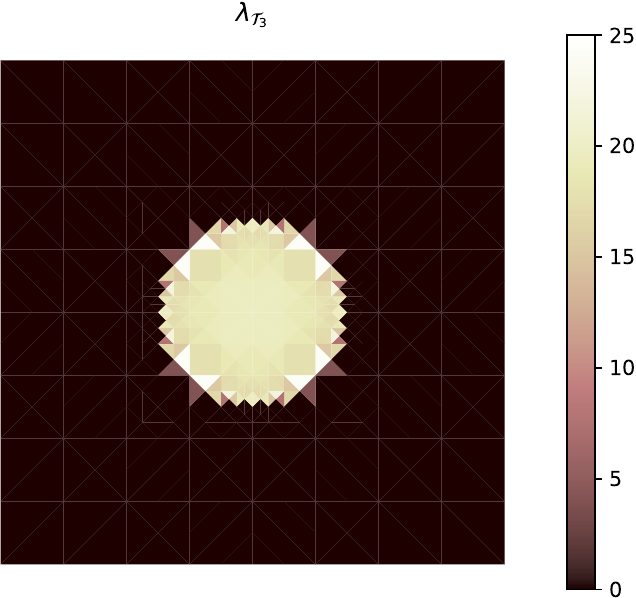}
  \includegraphics[width=0.33\textwidth]{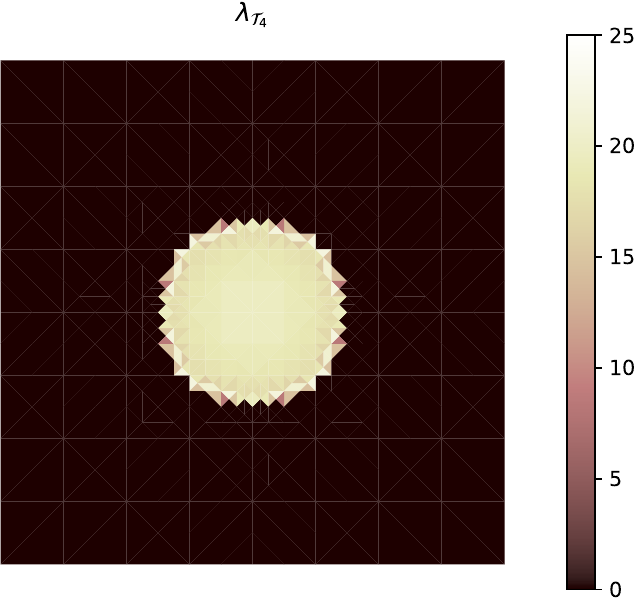}
  \includegraphics[width=0.33\textwidth]{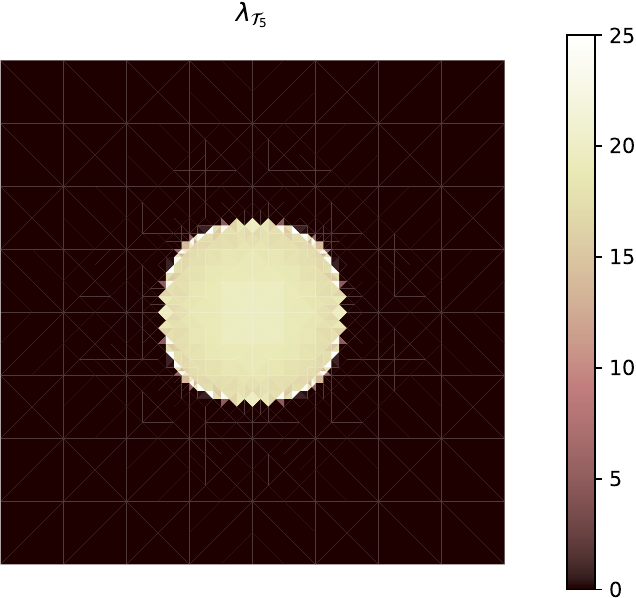}
\caption{The sequence of Lagrange multipliers for the membrane problem.}
\label{fig:membranelambdas}
\end{figure}

\begin{figure}
  \includegraphics[width=0.5\textwidth]{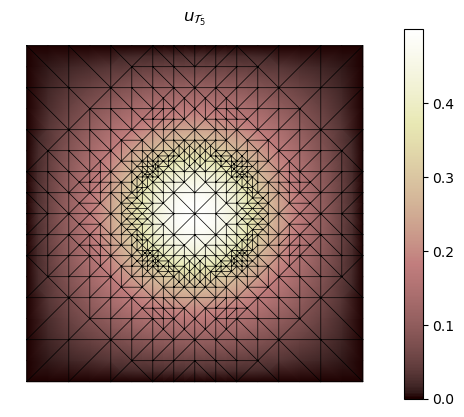}
\caption{The displacement of the membrane.}
\label{fig:membraneus}
\end{figure}

We will start with an academic example
of a membrane contact problem, i.e.~we solve
\eqref{eq:memobs} with $\Omega = (0, 1)\times (0,1)$,
$f=0$, and $G=1$. The obstacle is given by
\[
g = g(x,y) = \sin \pi x  \sin \pi y - \tfrac12.
\]
The resulting sequence of error indicators and adaptive meshes are
given in Figure~\ref{fig:membraneetas} and Figure~\ref{fig:membranemeshes}.
The solutions are given in Figure~\ref{fig:membranelambdas}
and Figure~\ref{fig:membraneus}.

We observe that the contact boundary is -- while not exactly a circle --
circular in shape which is a consequence of the symmetry of $g$
with respect to the axes $x=\tfrac12$ and $y=\tfrac12$.
The error indicators in Figure~\ref{fig:membraneetas}
show that depending on the iteration, different
parts of the error indicator are active and, hence,
responsible for the refinement of particular areas within the computational domain.
For instance, it becomes clear that $\eta_{C,K}$ is mainly
responsible for refinement within and on the boundary of
the coincidence set $\Omega_C$ while $\eta_{\partial K}$ is more
active outside of $\Omega_C$.

\begin{figure}
  \includegraphics[width=0.8\textwidth]{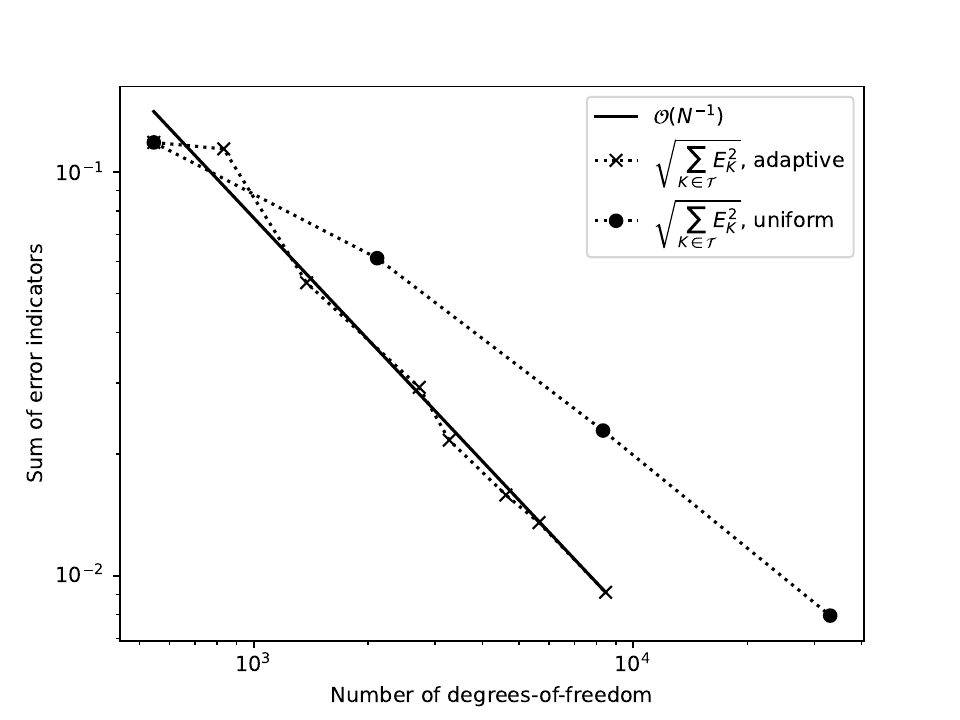}
\caption{The sum of the error indicators for the membrane problem.}
\label{fig:membraneerror}

\end{figure}

The shape of the coincidence set is clearly visible in the contact
force, Figure~\ref{fig:membranelambdas}, whose approximation is visually
improved by the use of adaptive meshing.
The discrete solution $u_{\Th_5}$ corresponding to the fifth
mesh is given in Figure~\ref{fig:membraneus}.
The best approach for observing the convergence of the
adaptive iteration is to plot the sum of the error indicators
as a function of the degrees-of-freedom $N$,
see Figure~\ref{fig:membraneerror} which also compares
the result between adaptive and uniform mesh sequences.
Based on the a posteriori error estimate, Theorem~\ref{thm:aposteriori},
the sum of error indicators gives an upper bound for the total error.
For the adaptive mesh sequence, we can expect the total error to
always converge at a rate which derives directly
from the polynomial degree of the finite element spaces
while the uniform mesh sequence remains limited by the
regularity of the loading, the boundary data or the
shape of the computational domain.

For quadratic finite element basis and two-dimensional smooth problems
the error is $\mathcal{O}(N^{-1})$.
For the uniform mesh sequence the error is $\mathcal{O}(N^{-\alpha})$
where $\alpha > 0$ may have any value less or equal
to $1$ depending on the problem setup
while the adaptive scheme
recovers the rate $\mathcal{O}(N^{-1})$ irrespective
of the smoothness of the loading, boundary data or the
shape of the boundary.
In addition to possible improvements in the asymptotic rates,
we usually observe a clear improvement in the absolute error.
For example, based on the rightmost data points in Figure~\ref{fig:membraneerror},
we can obtain the same discretisation error
with almost $75\,\%$ reduction in the number of degrees-of-freedom.

\begin{figure}
  \includegraphics[width=0.32\textwidth]{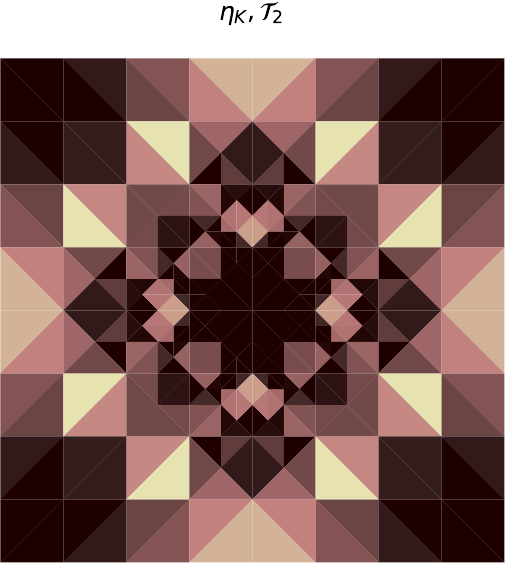}
  \includegraphics[width=0.32\textwidth]{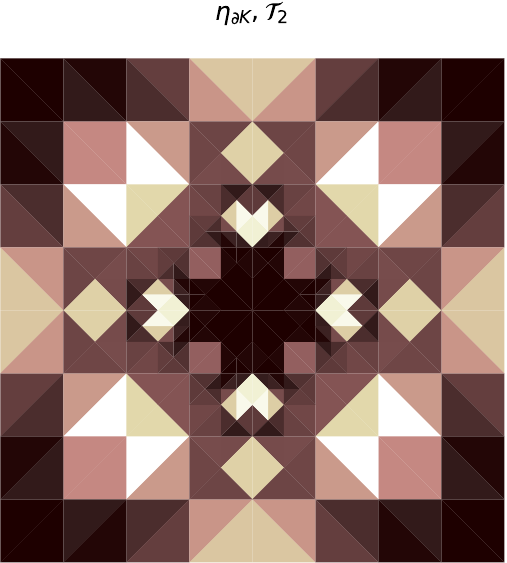}
  \includegraphics[width=0.32\textwidth]{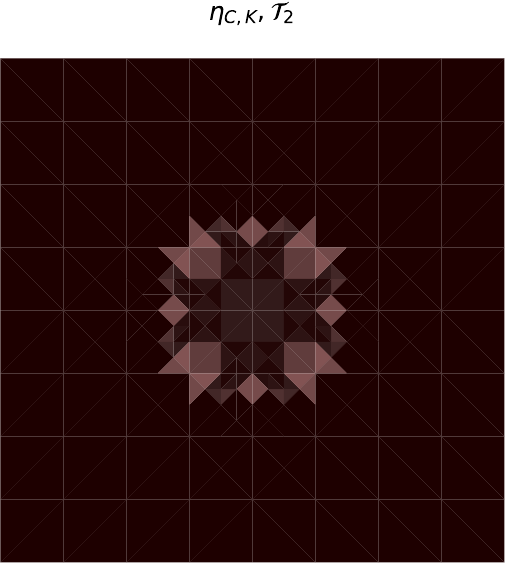}\\[0.2cm]
  \includegraphics[width=0.32\textwidth]{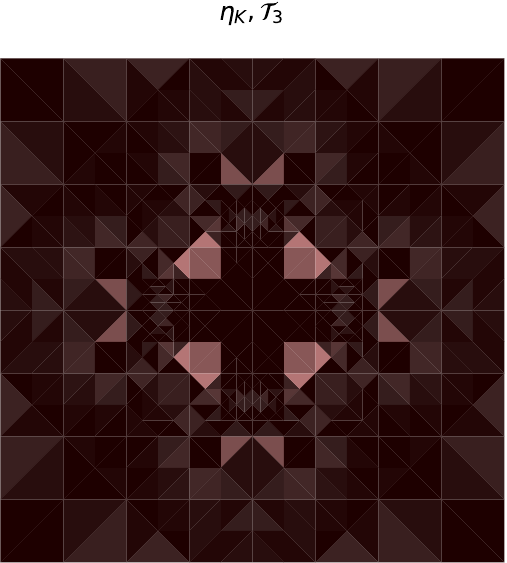}
  \includegraphics[width=0.32\textwidth]{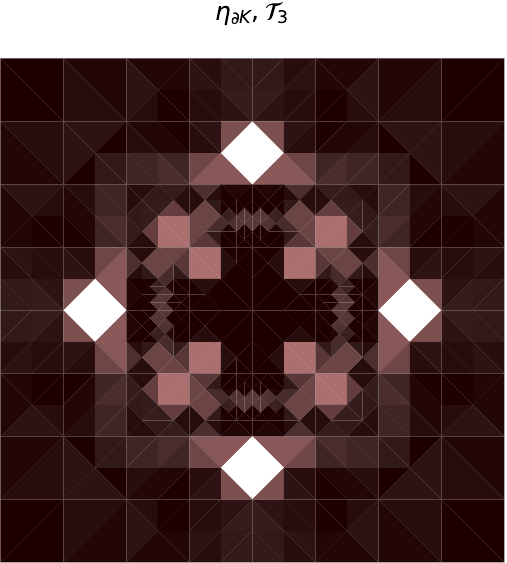}
  \includegraphics[width=0.32\textwidth]{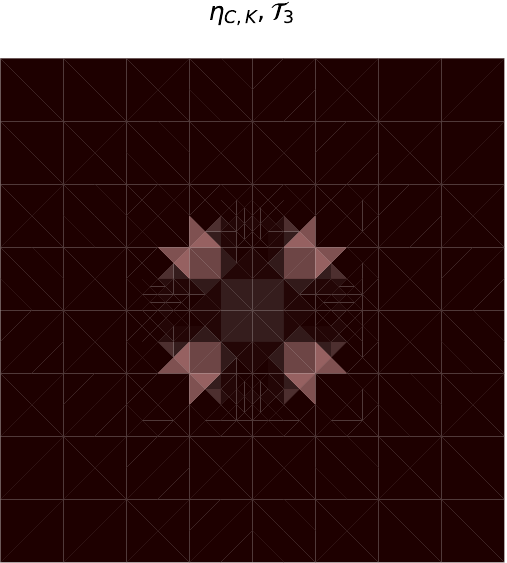}
\caption{(The sequence of error indicators for the membrane problem, part 2/3)}
\label{fig:membraneetas2}
\end{figure}

\begin{figure}
  \includegraphics[width=0.32\textwidth]{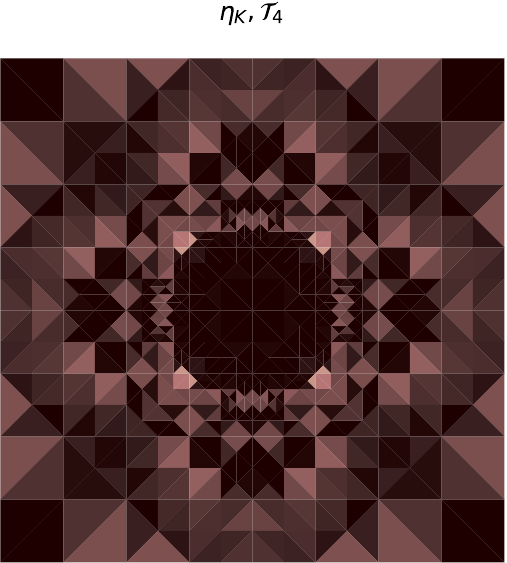}
  \includegraphics[width=0.32\textwidth]{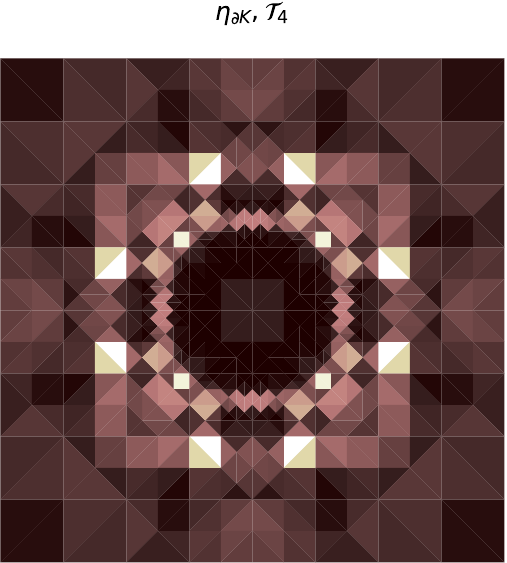}
  \includegraphics[width=0.32\textwidth]{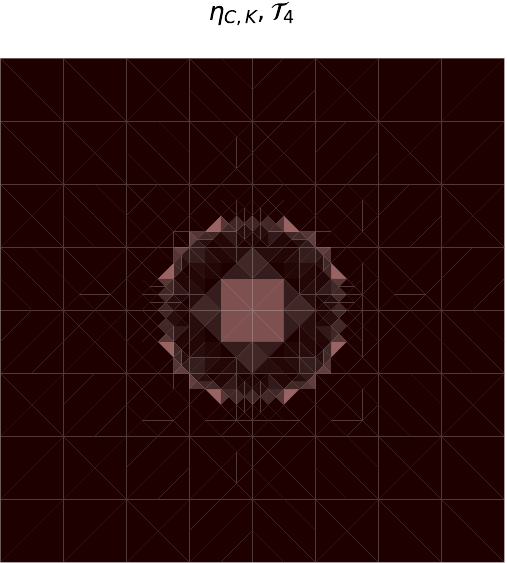}\\[0.2cm]
  \includegraphics[width=0.32\textwidth]{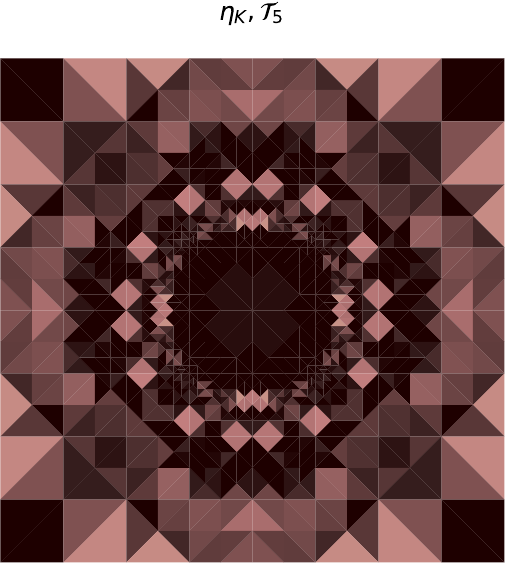}
  \includegraphics[width=0.32\textwidth]{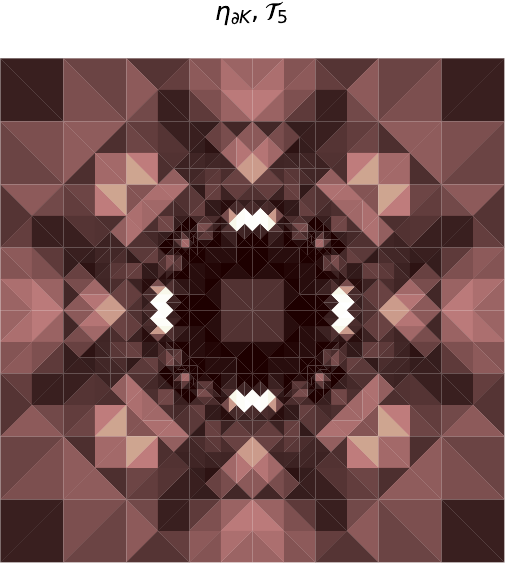}
  \includegraphics[width=0.32\textwidth]{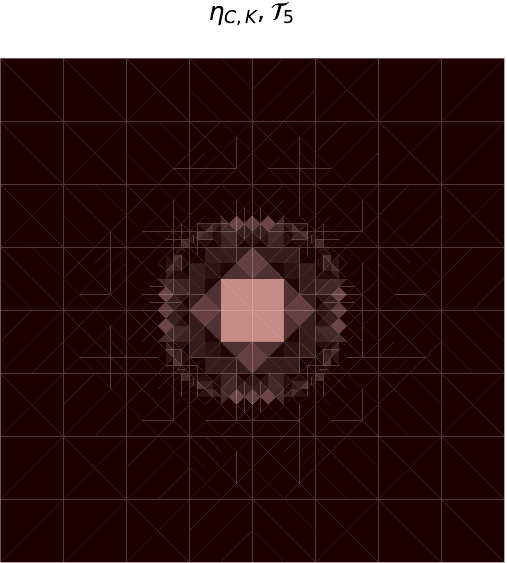}
\caption{(The sequence of error indicators for the membrane problem, part 3/3)}
\label{fig:membraneetas3}
\end{figure}

\newpage

\subsection{Plastic torsion of I beam}

\begin{figure}
  \includegraphics[width=0.5\textwidth]{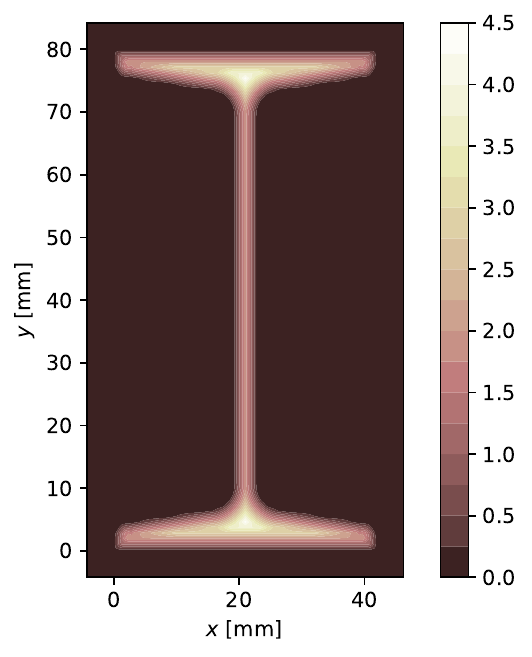}
  \includegraphics[width=0.35\textwidth]{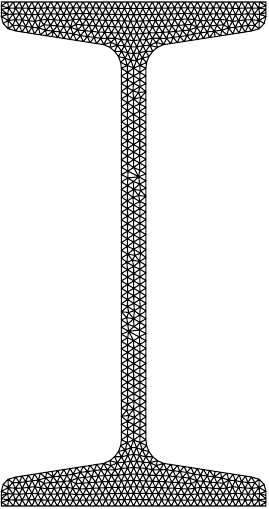}
\caption{The distance field for the IPN 80 profile and the initial mesh.}
\label{fig:torsiondf}
\end{figure}

Next we consider the torsion of a steel I beam
with standard European IPN 80 profile.
The cross section boundary is modelled using distance
fields which can be used both for mesh generation~\cite{persson2004simple}
and for defining $\delta$ in the formulation \eqref{eq:torobs}.
The distance field and the initial mesh are depicted in Figure~\ref{fig:torsiondf}.
The parameters of the problem are fixed in Table~\ref{tab:torsionparams} and
the twist per unit length $\theta$ is calculated as
\[
\theta = \frac{2\pi \gamma}{180^\circ l}.
\]
This problem is modelled after the experimental setup in Boru~\cite{boru2023}
which allows verifying that the plastic deformation occurs
roughly at the same angle of twist as in the experimental
setup.

\begin{table}
\caption{Parameters of the plastic torsion problem.}
\label{tab:torsionparams}%
\centering
\begin{tabular}{|l|l|}
\hline
\emph{parameter (symbol)} & \emph{value (unit)} \\ \hline
shear modulus ($G$) & 79.3 (GPa) \\ \hline
yield stress ($\tau$) & 0.240 (GPa) \\ \hline
angle of twist ($\gamma$) & 8 ($^\circ$) \\ \hline
length of the beam ($l$) & 400 (mm) \\ \hline
\end{tabular}
\end{table}

\begin{figure}
  \includegraphics[width=0.6\textwidth]{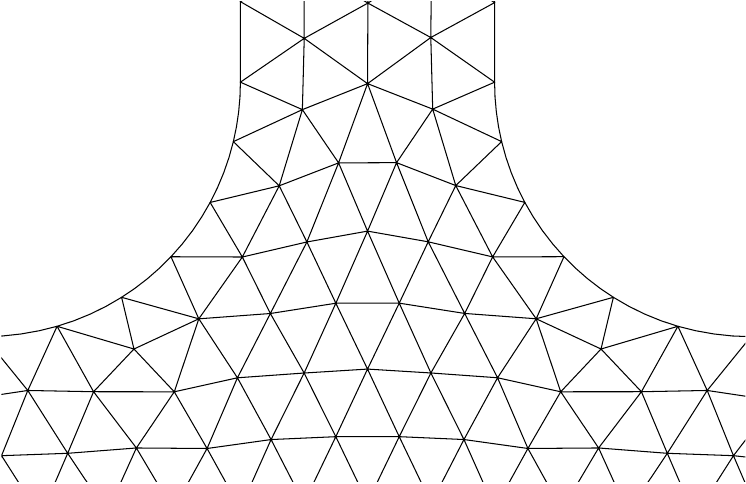}
\caption{A close-up view of the quadratic initial mesh.}
\label{fig:torsioncloseup}
\end{figure}

\begin{figure}
  \includegraphics[width=0.2\textwidth]{adaptive_torsion_0_mesh-crop.pdf}
  \includegraphics[width=0.2\textwidth]{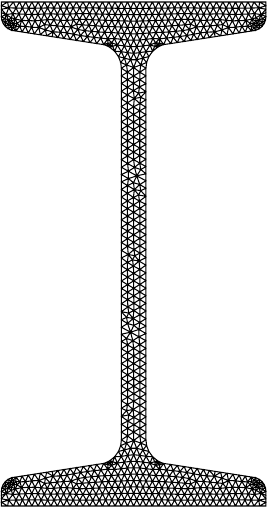}
  \includegraphics[width=0.2\textwidth]{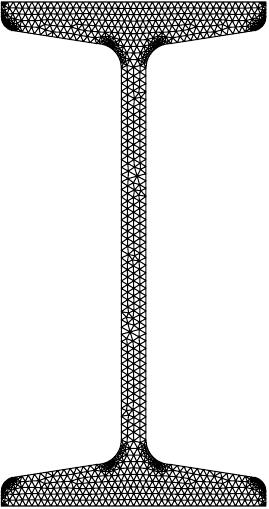}\\
  \includegraphics[width=0.2\textwidth]{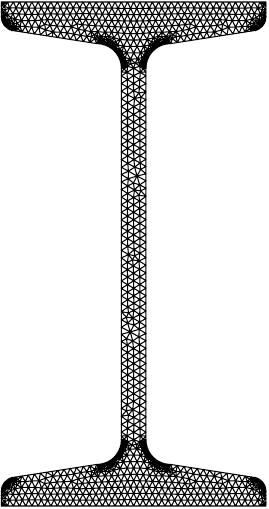}
  \includegraphics[width=0.2\textwidth]{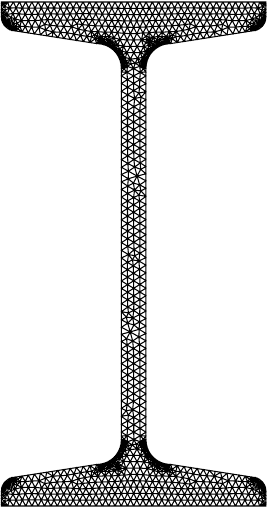}
  \includegraphics[width=0.2\textwidth]{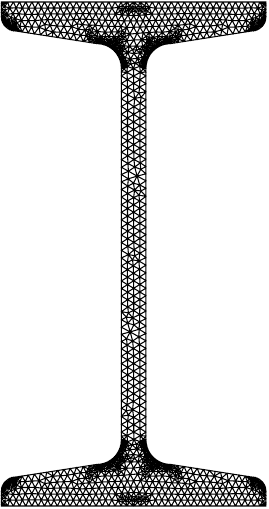}\\
  \includegraphics[width=0.2\textwidth]{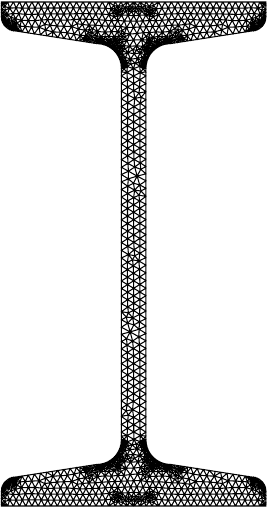}
  \includegraphics[width=0.2\textwidth]{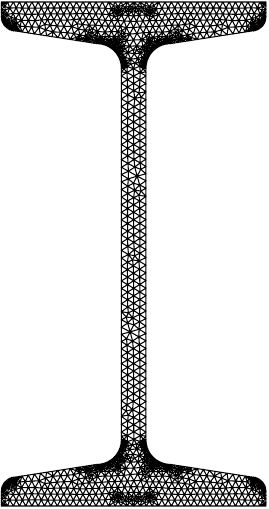}
  \includegraphics[width=0.2\textwidth]{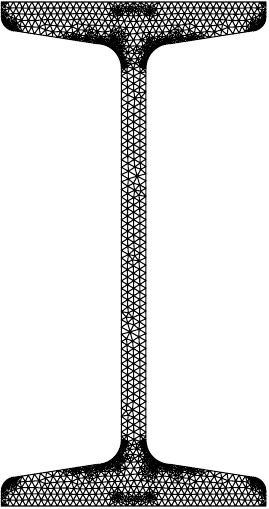}
\caption{The sequence of meshes for the torsion problem, part 1/2.}
\label{fig:torsionmeshes}
\end{figure}

\begin{figure}
  \includegraphics[width=0.2\textwidth]{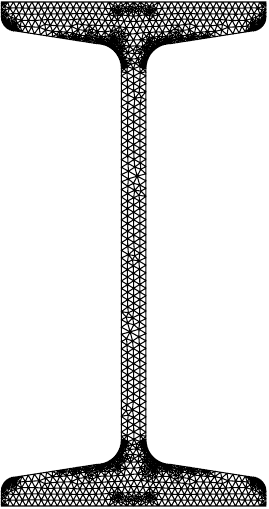}
  \includegraphics[width=0.2\textwidth]{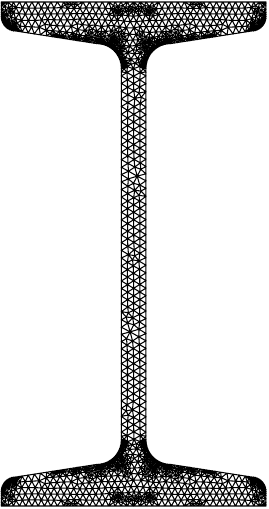}
  \includegraphics[width=0.2\textwidth]{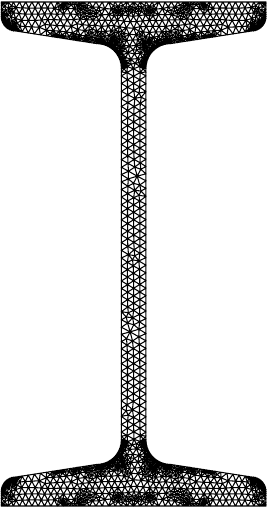}\\
  \includegraphics[width=0.2\textwidth]{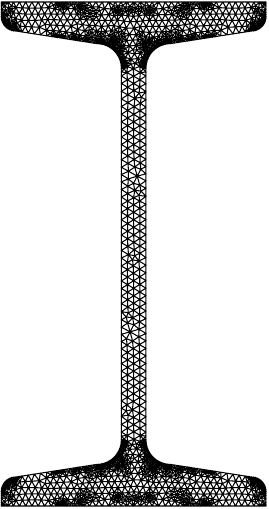}
  \includegraphics[width=0.2\textwidth]{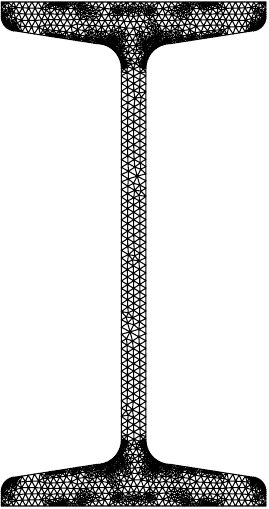}
  \includegraphics[width=0.2\textwidth]{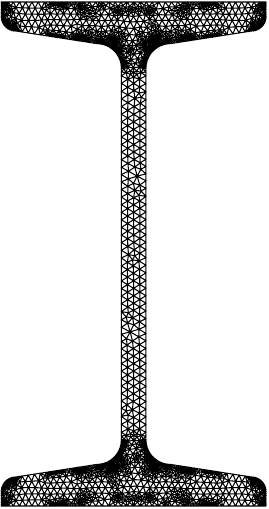}\\
  \includegraphics[width=0.2\textwidth]{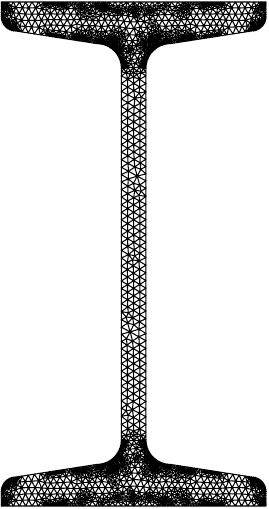}
  \includegraphics[width=0.2\textwidth]{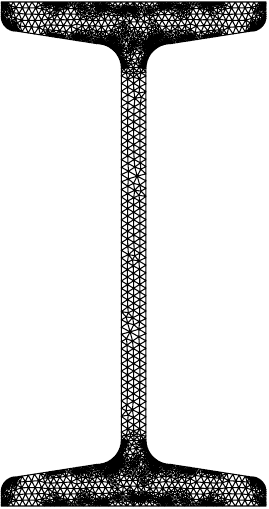}
  \includegraphics[width=0.2\textwidth]{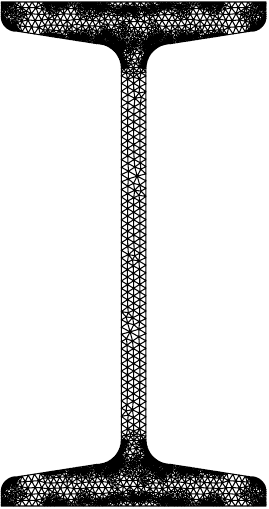}
\caption{The sequence of meshes for the torsion problem, part 2/2.}
\label{fig:torsionmeshes2}
\end{figure}

This is a challenging problem for the adaptive scheme because
parts of the domain boundary consist of circle arcs.
Consequently, we have implemented the following changes
to the basic adaptive scheme:
\begin{enumerate}
\item The mesh uses quadratic isoparametric mapping to approximate
  the domain boundary using piecewise quadratic segments (see Figure~\ref{fig:torsioncloseup})
\item The mesh refinement algorithm uses the distance field of Figure~\ref{fig:torsiondf}
  to project any new mesh nodes onto the boundary
\item Laplacian mesh smoothing is applied after each refinement step.
\end{enumerate}
The purpose of the change 1 is to minimise the effect of the variational
crime due to inexact boundary representation.
It is not possible to completely remove this effect using a piecewise polynomial mapping
because circles cannot be represented exactly by polynomials --
the consequence of which is exemplified by the sequence
of meshes in Figure~\ref{fig:torsionmeshes} and Figure~\ref{fig:torsionmeshes2}
with some obvious refinement also near the
end of the flanges consisting of outer circular arcs.
Fortunately, the adaptive scheme is clearly not stuck on refining a single
part of the domain only which would be the case if linear boundary
segments were used.

The change 2 is also essential and might be challenging to implement
for arbitrary domains and meshes.  However, as our domain 
is already defined using the distance field $\delta$, it becomes a simple
evaluation of the distance field at the locations of the
new nodes and moving them in the direction of the normal vector.
The change 3 is optional but retaining mesh regularity
in the red-green-blue refinement
requires, in general, checking how many times each edge
has been split and changing
the refinement type from green to blue accordingly.
This may in some cases lead to
a blow-up of the refinement zone size~\cite[p.~68]{verfurthbook}
and, hence, we avoid changing the refinement
type and instead use Laplacian smoothing to improve
mesh regularity.

\begin{figure}
  \includegraphics[width=0.8\textwidth]{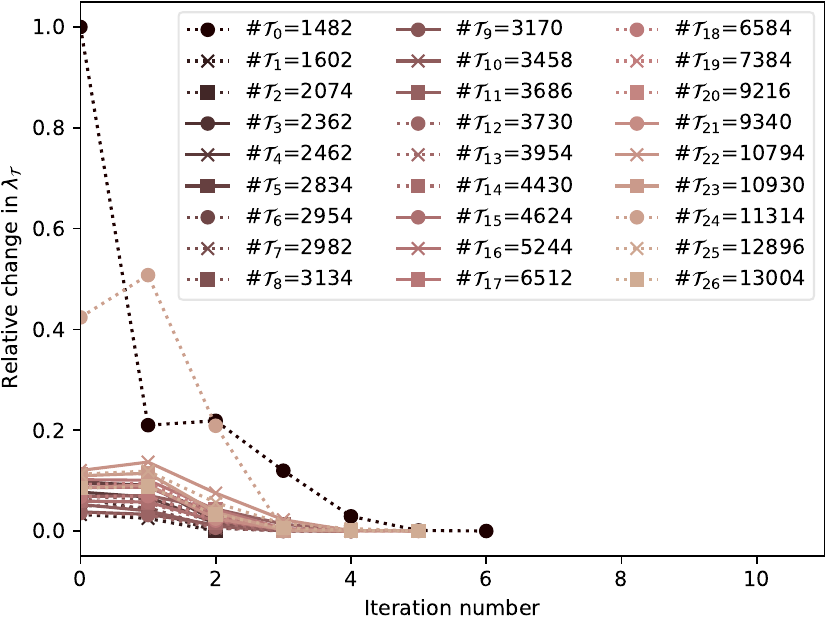}
\caption{Primal-dual active set method.}
\label{fig:torsionprimaldual}
\end{figure}

\begin{figure}
  \includegraphics[width=0.25\textwidth]{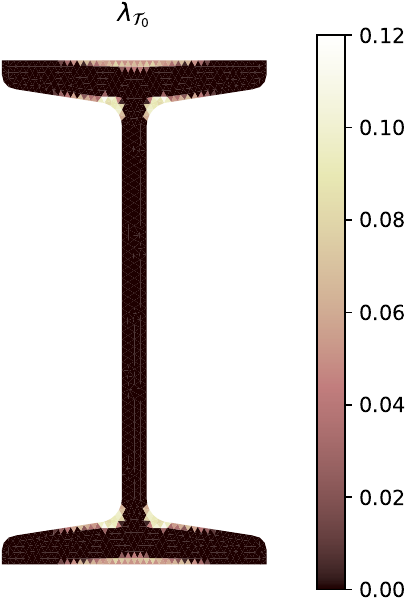}
  \includegraphics[width=0.25\textwidth]{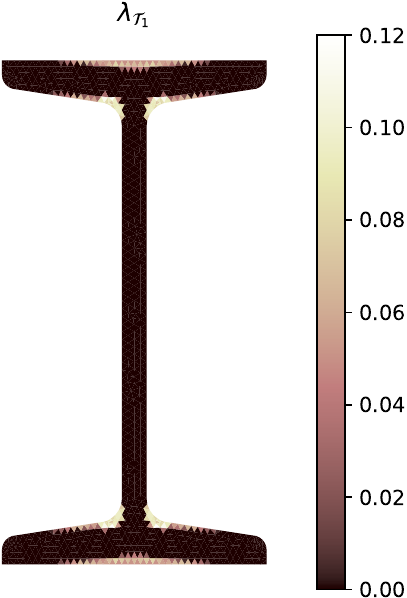}
  \includegraphics[width=0.25\textwidth]{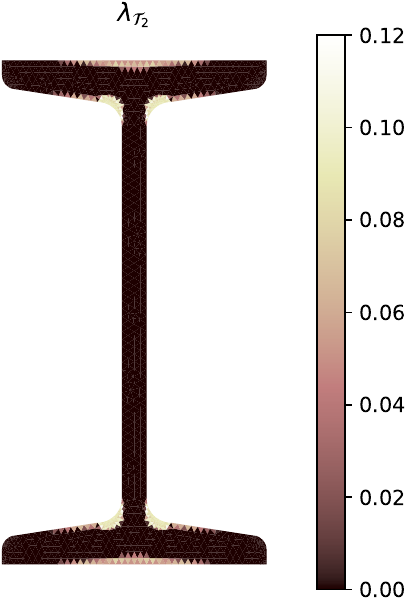}\\
  \includegraphics[width=0.25\textwidth]{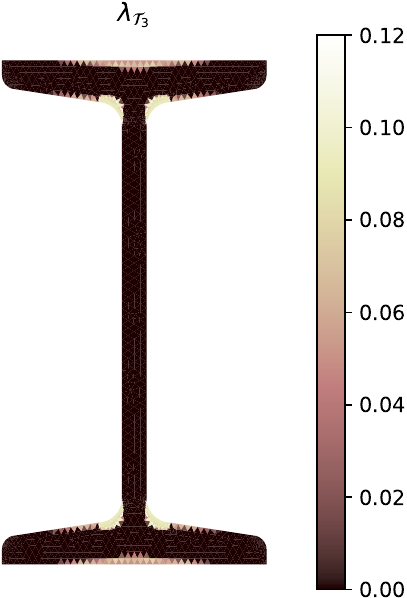}
  \includegraphics[width=0.25\textwidth]{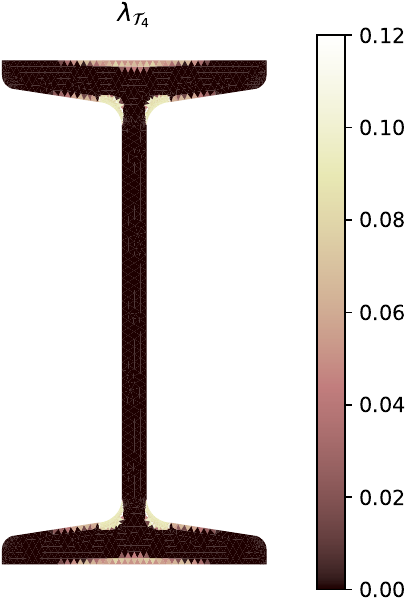}
  \includegraphics[width=0.25\textwidth]{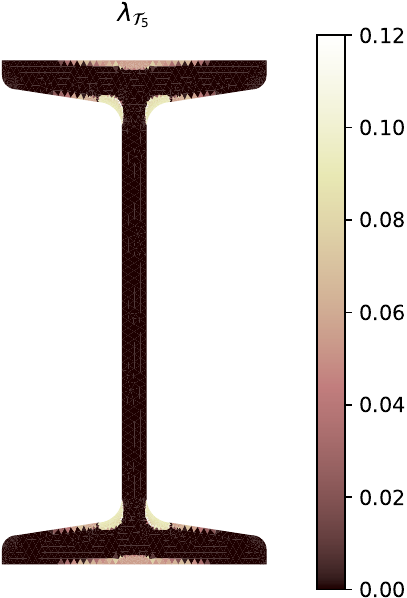}\\
  \includegraphics[width=0.25\textwidth]{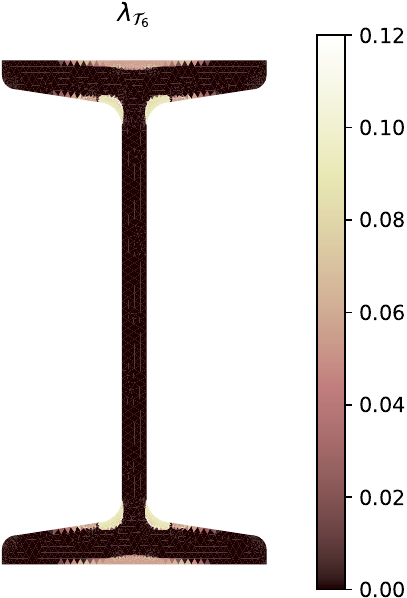}
  \includegraphics[width=0.25\textwidth]{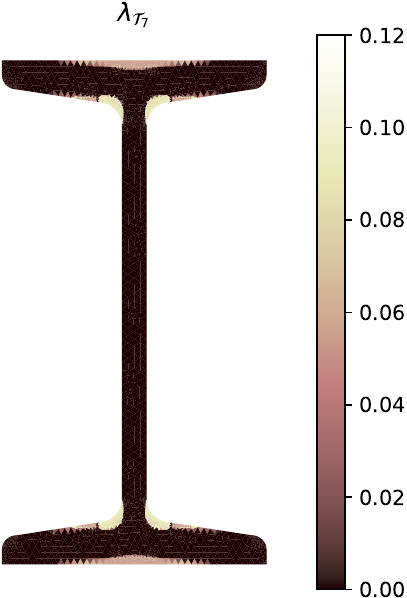}
  \includegraphics[width=0.25\textwidth]{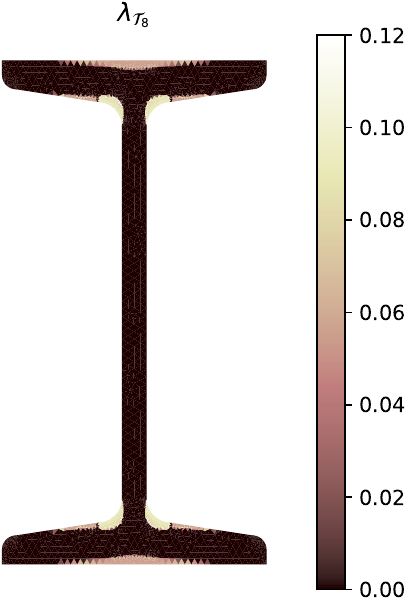}
\caption{The sequence of Lagrange multipliers for the torsion problem, part 1/2.}
\label{fig:torsionlambdas}
\end{figure}

\begin{figure}
  \includegraphics[width=0.25\textwidth]{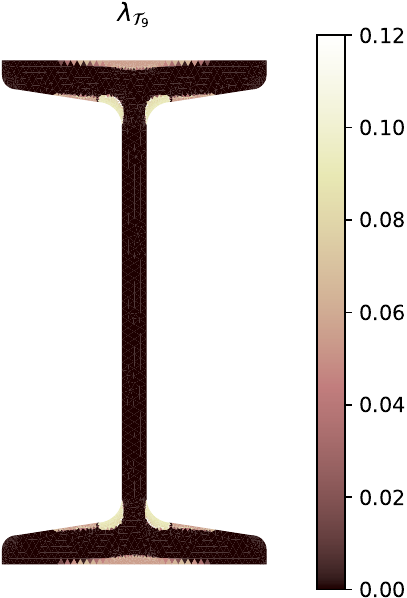}
  \includegraphics[width=0.25\textwidth]{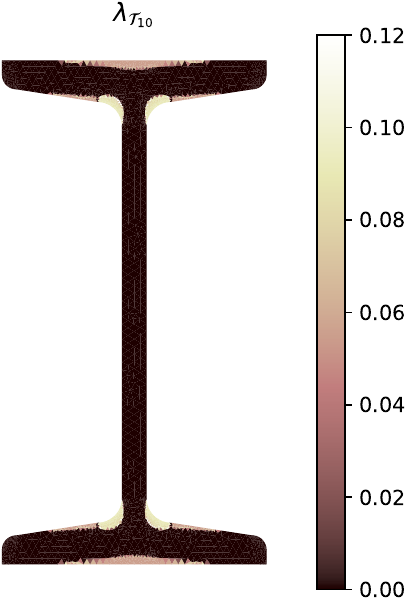}
  \includegraphics[width=0.25\textwidth]{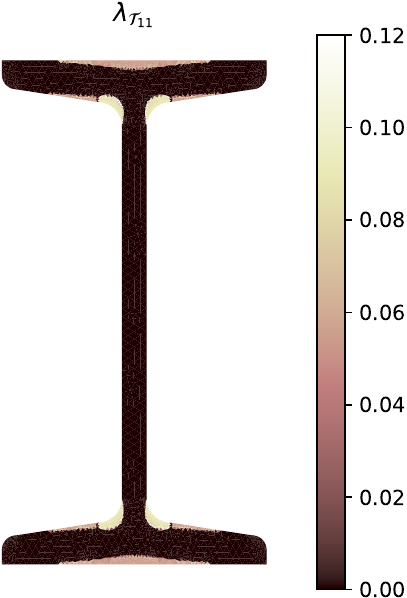}\\
  \includegraphics[width=0.25\textwidth]{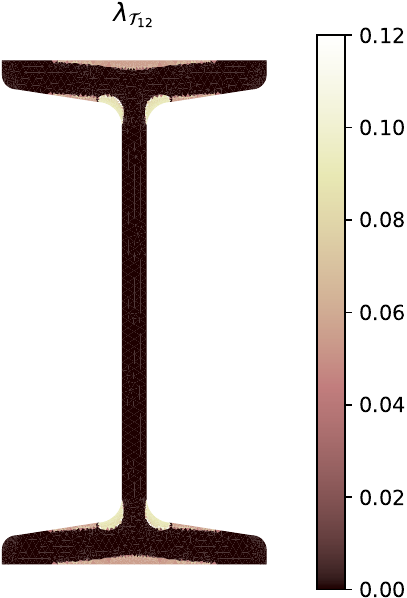}
  \includegraphics[width=0.25\textwidth]{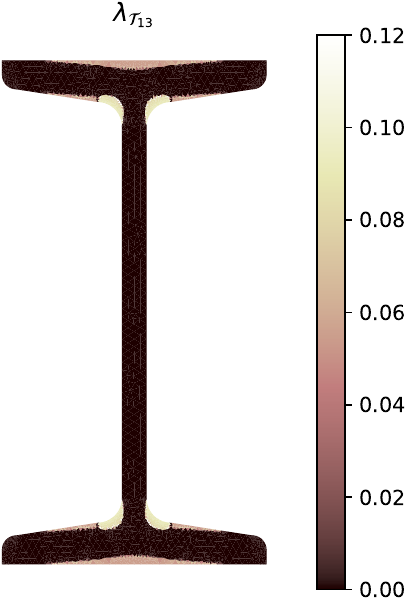}
  \includegraphics[width=0.25\textwidth]{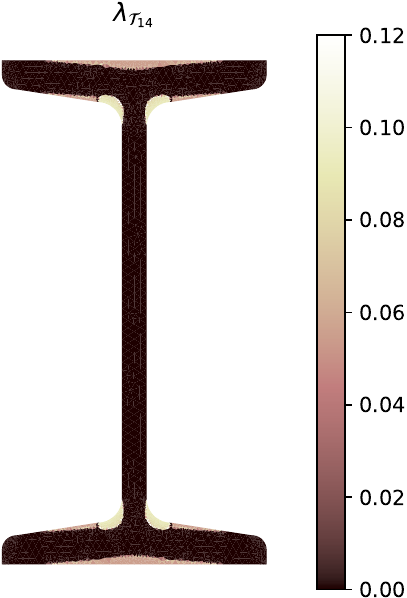}\\
  \includegraphics[width=0.25\textwidth]{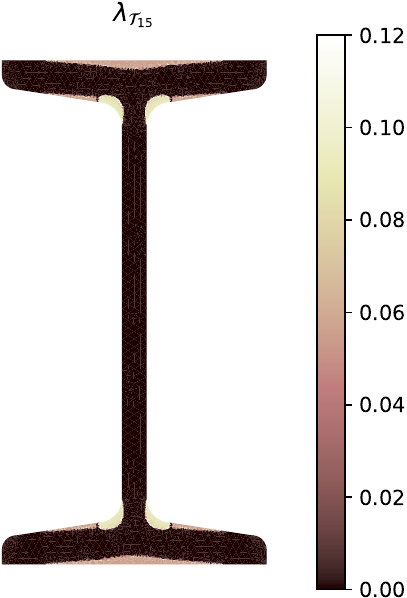}
  \includegraphics[width=0.25\textwidth]{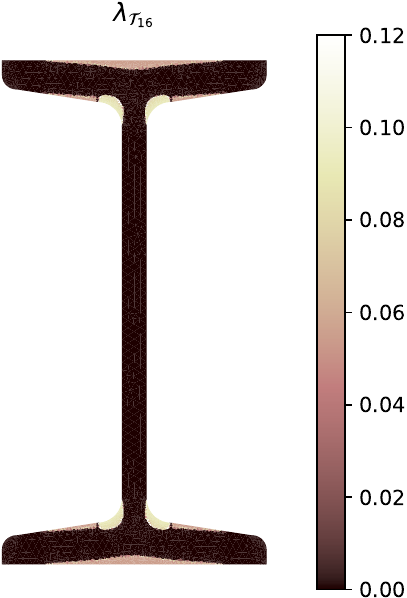}
  \includegraphics[width=0.25\textwidth]{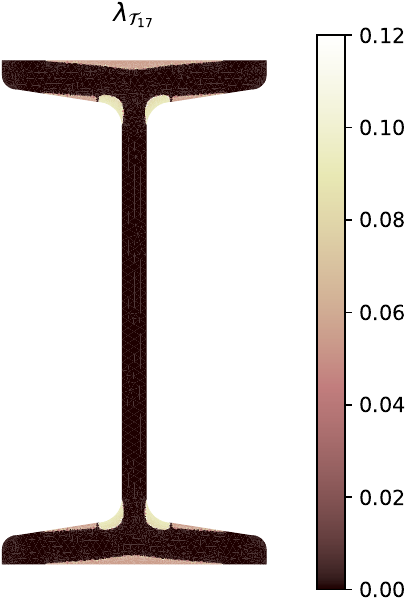}
\caption{The sequence of Lagrange multipliers for the torsion problem, part 2/2.}
\label{fig:torsionlambdas2}
\end{figure}

\begin{figure}
  \includegraphics[width=0.7\textwidth]{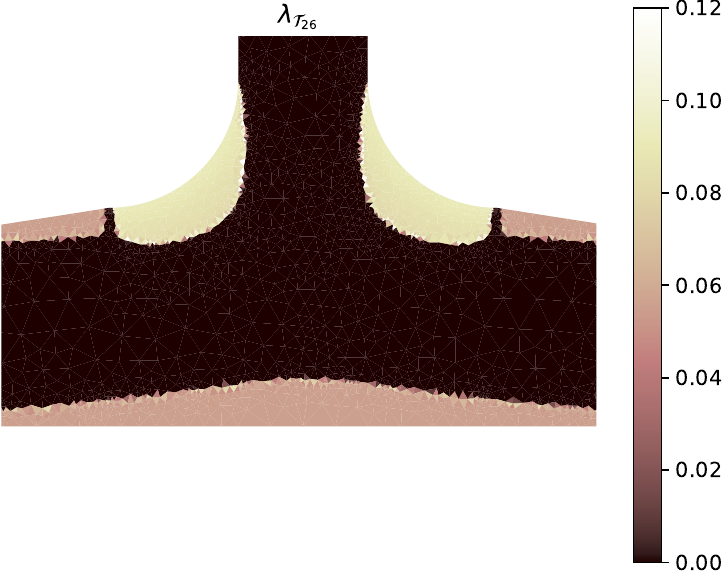}
\caption{A close-up view of the Lagrange multiplier.}
\label{fig:torsionzones}
\end{figure}

\begin{figure}
  \includegraphics[width=0.7\textwidth]{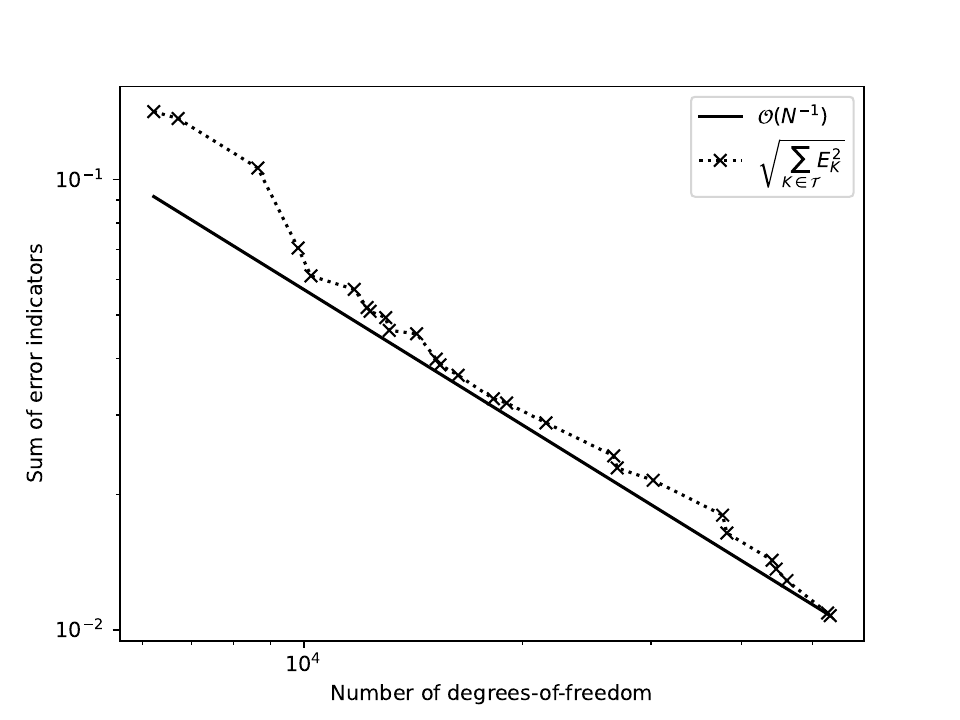}
\caption{The sum of the error indicators for the torsion problem.}
\label{fig:torsionerror}
\end{figure}

We note that the use of mesh smoothing will make it difficult and
computationally more demanding to reuse the solution as an initial guess.
In particular, the projection between nonmatching meshes
requires finding intersections of the triangles for constructing
integration points and calculating $L^2$ projection from
one mesh to another.
After this implementational challenge is solved,
we are able to again bound
the number of primal-dual interations by a small constant, see Figure~\ref{fig:torsionprimaldual}.
The first 18 solutions from the sequence of Lagrange multipliers are
given in Figure~\ref{fig:torsionlambdas} and Figure~\ref{fig:torsionlambdas2},
and a close-up view of the plastic zones corresponding
to the final solution is given in Figure~\ref{fig:torsionzones}.
While the initial mesh does not
distinguish between the different plastic zones,
the different zones are clearly visible
in the final solution.
The convergence is again tracked using
the sum of the error indicators, given in Figure~\ref{fig:torsionerror}.

\newpage

\subsection{Cavitation in a hydrodynamic bearing}

As a final computational demonstration, we consider
the lubricant pressure in a hydrodynamic bearing.
The problem is modelled after
the experimental setup presented in
Etsion--Ludwig~\cite{etsion1982} and
applied further in
Braun--Hendricks~\cite{braun1984},
see Table~\ref{tab:bearingparams}
for the parameters we have used
in the simulation.
We note that
the lubricant viscosity
is highly dependent on the temperature
of the lubricant
and that the temperature is not constant
within the bearing \cite{etsion1982}.
However, based on the measurement data,
it becomes obvious that the
bearing rotational speed and the
lubricant temperature have a correlation \cite{braun1984}.
This can be explained by the fact that higher
rotational speeds lead to higher shear rates
and, hence, higher viscous heating
within the lubricant.

\begin{figure}
  \includegraphics[width=0.7\textwidth]{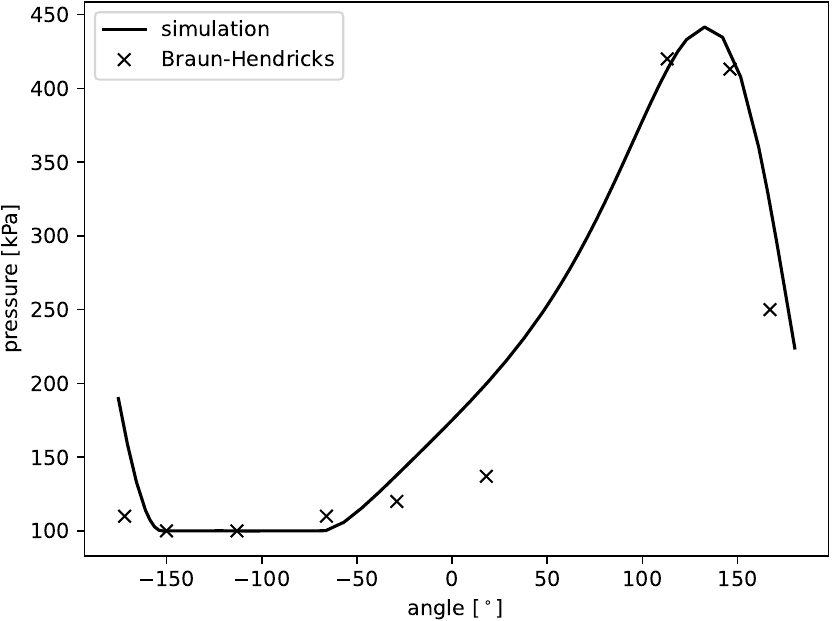}
\caption{(The pressure at the midline) A comparison is done against a few experimental data points from Braun--Henricks~\cite{braun1984} which contains the measured pressure distribution over the entire bearing surface.}
\label{fig:bearingcomparison}
\end{figure}

Consequently, the
experimental maximum pressure \cite{braun1984} is
not linearly proportional to the surface velocity $V$,
as one would expect based
on the right hand side of \eqref{eq:reynolds},
since the temperature
is inversely proportional to the
viscosity of the lubricant $\mu$.
There exist separate models for thermohydrodynamics
that could be utilised for studying the effect of the
nonzero temperature gradient;
see,
e.g., Feng at al.~\cite{feng2019}.
On the other hand, we could technically
use the temperature measurement data
presented in Braun--Hendricks~\cite{braun1984}
to calculate a space-dependent value
for the dynamic viscosity and use that
value in our simulation.
For simplicity, we choose a representative
constant value for $\mu$ which is close to the true
viscosity at the measured maximum temperature $33^\circ$C
while the measured minimum temperature
is approximately $22^\circ$C.
The cavitation pressure is a rough approximation
of the minimum pressure
in the pressure diagrams of \cite{braun1984}.
We note that this problem setup can
reasonably reproduce the experimental pressure at the midline, see
Figure~\ref{fig:bearingcomparison}
for a comparison.

\begin{figure}
  \includegraphics[width=0.49\textwidth]{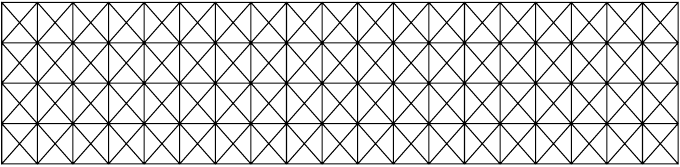}
  \includegraphics[width=0.49\textwidth]{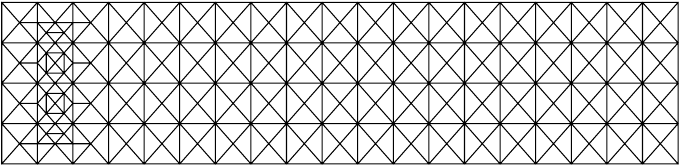}
  \includegraphics[width=0.49\textwidth]{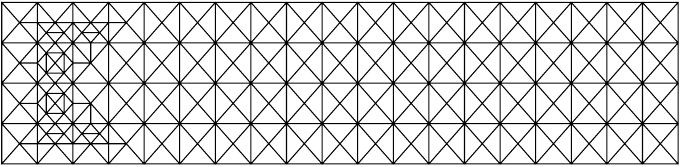}
  \includegraphics[width=0.49\textwidth]{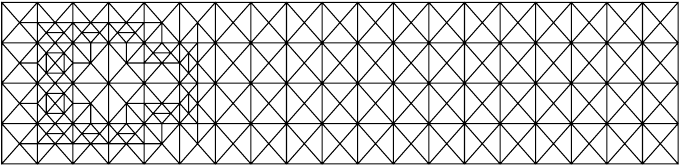}
  \includegraphics[width=0.49\textwidth]{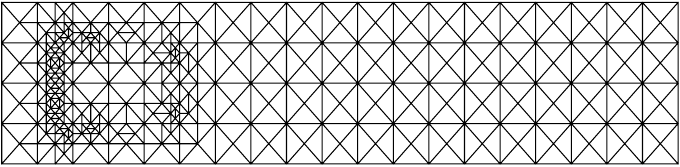}
  \includegraphics[width=0.49\textwidth]{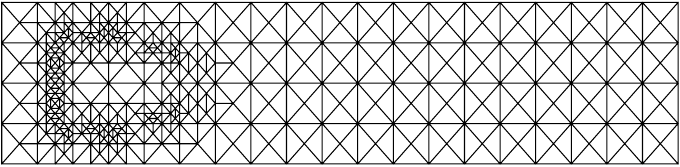}
  \includegraphics[width=0.49\textwidth]{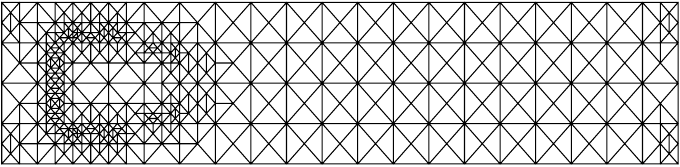}
  \includegraphics[width=0.49\textwidth]{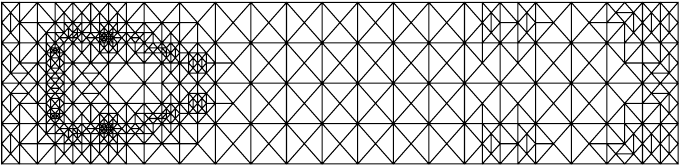}
  \includegraphics[width=0.49\textwidth]{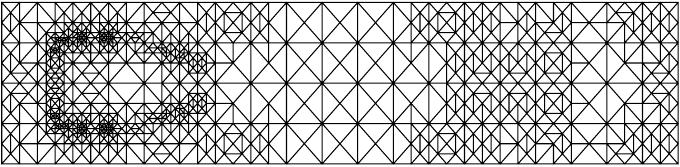}
\caption{The sequence of meshes for the bearing problem.}
\label{fig:bearingmeshes}
\end{figure}

\begin{figure}
  \includegraphics[width=0.49\textwidth]{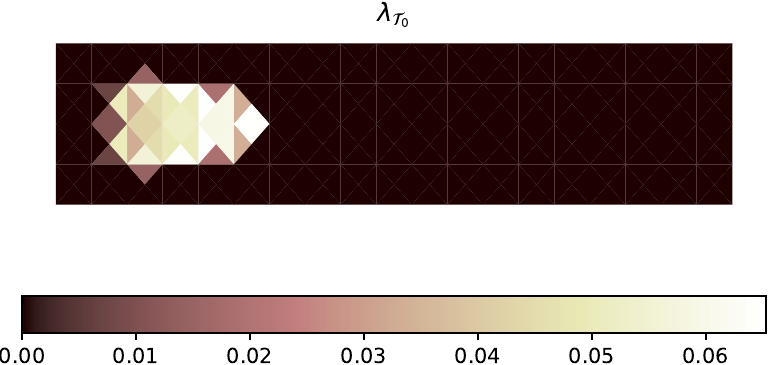}
  \includegraphics[width=0.49\textwidth]{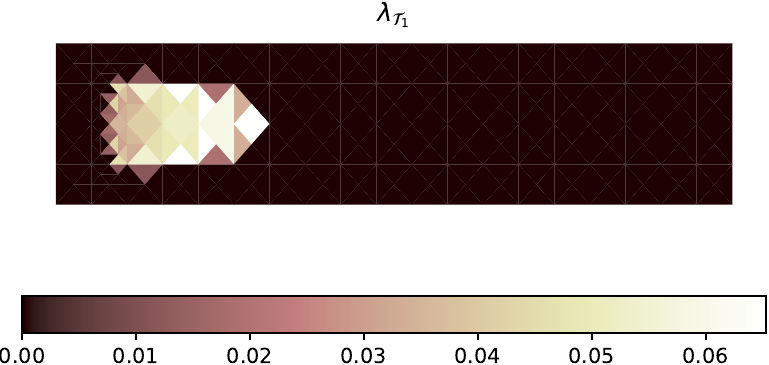}
  \includegraphics[width=0.49\textwidth]{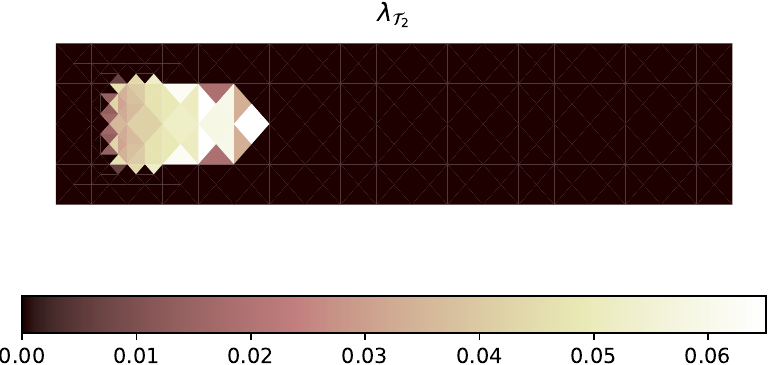}
  \includegraphics[width=0.49\textwidth]{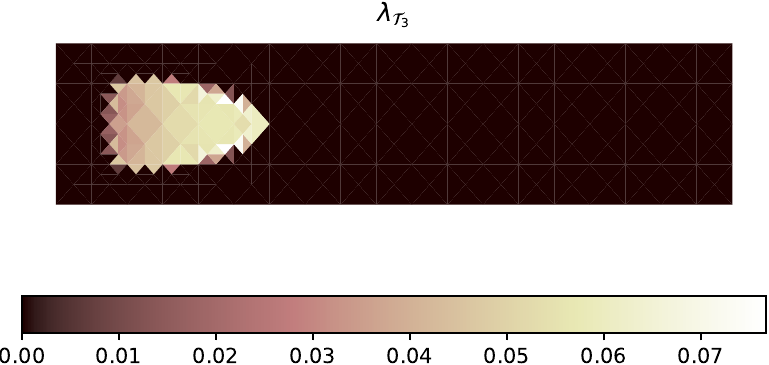}
  \includegraphics[width=0.49\textwidth]{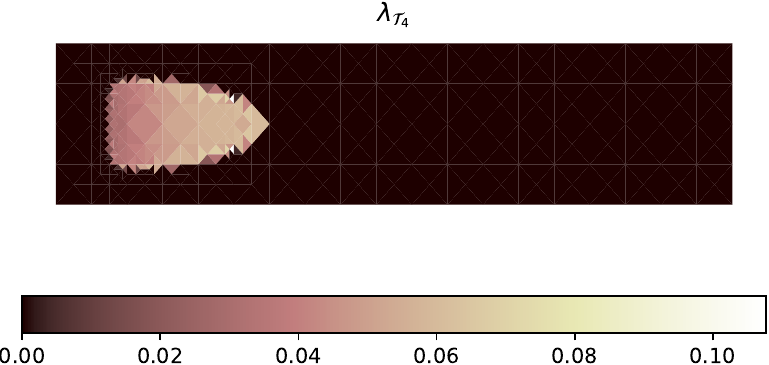}
  \includegraphics[width=0.49\textwidth]{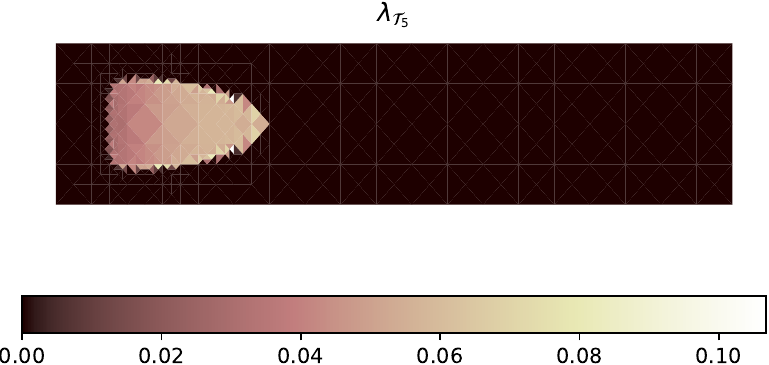}
  \includegraphics[width=0.49\textwidth]{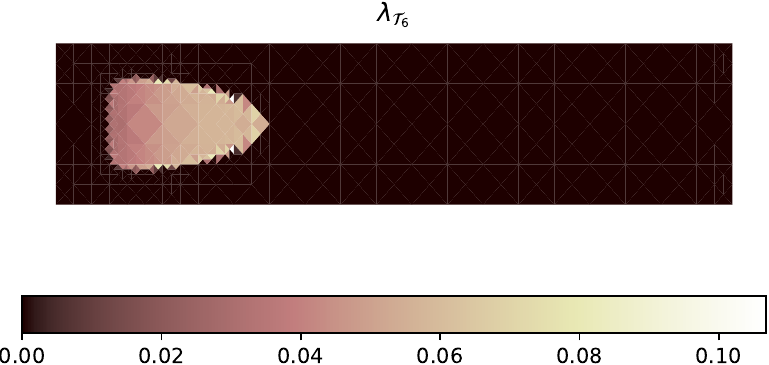}
  \includegraphics[width=0.49\textwidth]{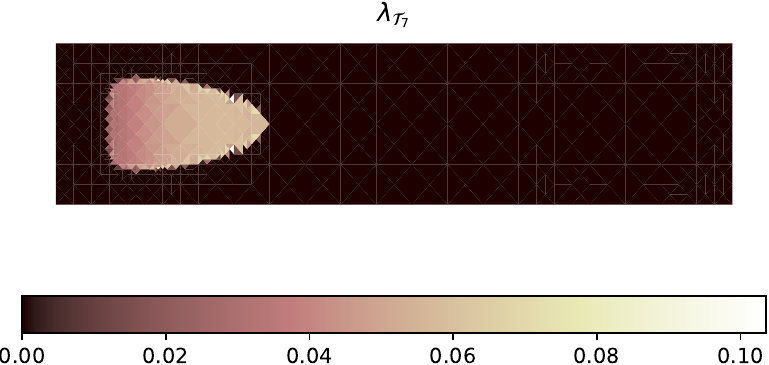}
  \includegraphics[width=0.49\textwidth]{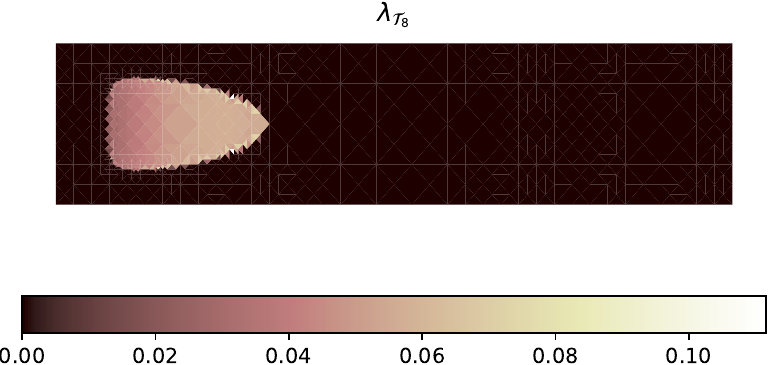}
\caption{The sequence of Lagrange multipliers for the bearing problem.}
\label{fig:bearinglambdas}
\end{figure}

\begin{figure}
  \includegraphics[width=0.6\textwidth]{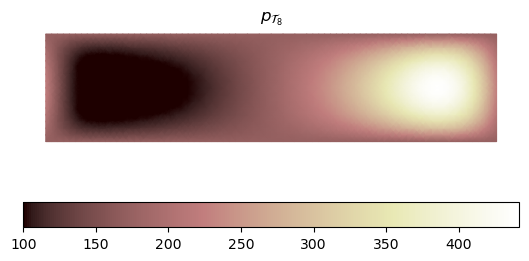}
\caption{The pressure for the bearing problem.}
\label{fig:bearingsols}
\end{figure}

\begin{figure}
  \includegraphics[width=0.7\textwidth]{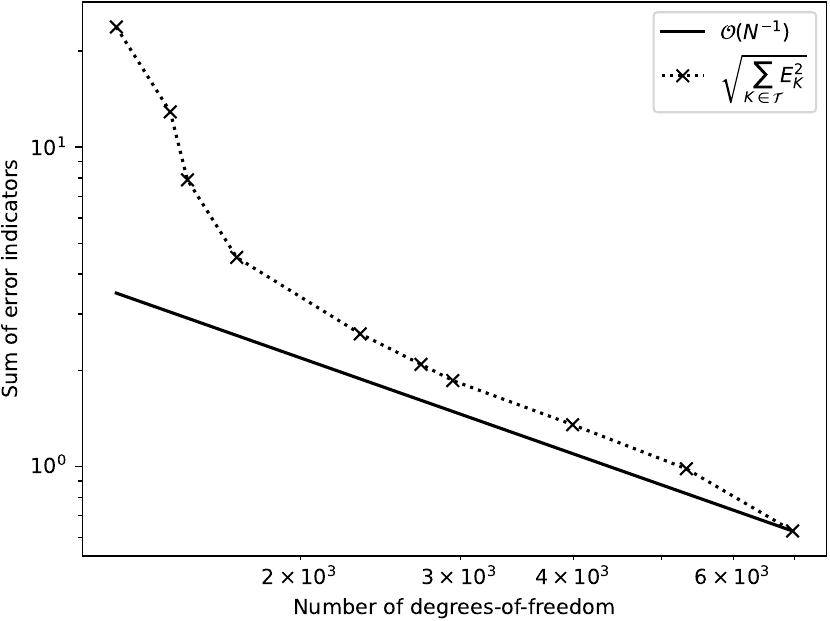}
\caption{The total error for the bearing problem.}
\label{fig:bearingerror}
\end{figure}

We do the following changes to the error indicator
to account for the large variations
of $d^3$ within the computational domain \cite{gustafsson2018}:
\begin{align*}
  \eta_K^2 &= \frac{h_K^2}{d_K^3} \| \nabla \cdot(\tfrac{d^3}{\mu} \nabla p_\Th) + \lambda_\Th + f \|_{L^2(K)}^2 \\
  \eta_{\partial K}^2 &= \frac 12 \frac{h_K}{d_K^3} \| \tfrac{d^3}{\mu} \llbracket \nabla p_\Th \cdot n \rrbracket\|_{L^2(\partial K \setminus \partial \Omega)}^2,
\end{align*}
where $d_K$ is the mean value of $d$ within the triangle $K$.
The sequence of adaptive meshes in Figure~\ref{fig:bearingmeshes}
shows that it is important, in terms of discretisation error,
to triangulate accurately the boundary of the cavitation zone $\Omega_C$.
We can observe the shape of the cavitation zone
in Figure~\ref{fig:bearinglambdas} and it
has a distincive shape, reminiscent of the photographs
taken of the experimental setup with transparent
bearing surfaces~\cite{braun1984}.
The pressure over the full journal bearing
is given in Figure~\ref{fig:bearingsols}
and the observed numerical convergence rate is depicted in
Figure~\ref{fig:bearingerror}.

\begin{table}
\caption{Parameters of the hydrodynamic bearing problem.}
\label{tab:bearingparams}%
\centering
\begin{tabular}{|l|l|l|}
\hline
\emph{parameter (symbol)} & \emph{value (unit)} & \emph{explanation} \\ \hline
radius of the bearing ($R$) & 25.4 (mm) & from \cite{braun1984} \\ \hline
length of the bearing ($L$) & 38.1 (mm) & from \cite{braun1984} \\ \hline
bearing clearance ($c_1$) & 114 ($\mu$m) & from \cite{etsion1982}  \\ \hline
dynamic viscosity at $T=20^\circ\text{C}$ ($\mu_{20^\circ\text{C}}$) & 0.0456 (Pa$\cdot$s) & from \cite{braun1984} \\ \hline
dynamic viscosity at $T=40^\circ\text{C}$ ($\mu_{40^\circ\text{C}}$) & 0.0232 (Pa$\cdot$s) & from \cite{braun1984} \\ \hline
dynamic viscosity ($\mu$) & 0.0307 (Pa$\cdot$s) & $\mu_{20^\circ\text{C}} < \mu < \mu_{40^\circ\text{C}}$  \\ \hline
eccentricity ($e$) & 0.4 (1) & from \cite{braun1984} \\ \hline
oil bath pressure ($p_{\text{env}}$) & 172 (kPa) & from \cite{braun1984} \\ \hline
cavitation pressure ($p_{\text{cav}}$) & 100 (kPa) & pressure diagrams in \cite{braun1984} \\ \hline
speed ($V$) & 5.320 (m/s) & $\approx$ 2000 RPM \cite{braun1984} \\ \hline
\end{tabular}
\end{table}

\newpage

\section{Discussion}

Adaptive finite elements use less degrees-of-freedom
to achieve a given level of accuracy.
Hence, they
allow increasing the limit of complexity
before supercomputing is required while
simultaneously reducing the user time spent on meshing.
This is also true for inequality-constrained problems,
such as the obstacle problems presented in this chapter,
or even more general elastic contact problems~\cite{gustafsson2020}.
For example, it is clear that 
the location of the coincidence set
may change in time-dependent problems and no fixed mesh can be ideal
over the entire time range.

On the other hand, there are general challenges related to the
practical adoption of adaptive methods including
but not limited to:
(1) the estimators are specific to a particular
continuous problem
and its discretisation, and different problems require different estimators
which means that the continuous formulation and its
relation to the discretisation must be well understood;
(2) nonpolygonal domains can be challenging or impossible to
represent exactly using standard finite elements
while a description of the exact geometry
must be available to the refinement algorithm;
and (3) different norms may be preferable, e.g.,
due to large variations in the material parameters or loadings,
or simply due to the preferences of the user,
that again lead to different error estimators.

The computational demonstrations presented in this chapter
exemplify some of the aforementioned challenges.
Many popular finite element solvers for day-to-day engineering
use fixed meshes and do not integrate geometry description,
as would be required to overcome challenge (2).
Solving challenges (1) and (3) in an automatic manner
would require an extensive catalogue of
estimators for various governing equations,
boundary conditions, material fields, and finite element methods.
For other problems of practical interest, such
as curved boundary, two body, or large displacement problems
in frictional contact~\cite{renard2013, chouly2014, gustafsson2022tresca, rodolfo2023},
the a~posteriori error analysis may not be developed
enough so that reliable adaptive schemes
can be established.

For the obstacle problems presented in this chapter,
the finite element error analysis is by now well understood
and adaptive methods can be used.
Presently, practical problems are most conveniently
implemented using programmable solvers that have
a built-in support for
adaptive mesh refinement with
quadratic or higher order geometry
description.
The examples in this chapter
were computed using scikit-fem~\cite{gustafsson2020scikit}
which is a library of assembly routines for implementing
finite element solvers in Python.
The mesh generation was done
using a variant of distmesh~\cite{persson2004simple}.
Other open source packages suitable
for implementing the methods
presented in this chapter include
FreeFem++~\cite{hecht2012},
NGSolve~\cite{gangl2021}
and FEniCSx~\cite{baratta2023}.

\bibliographystyle{unsrt}
\bibliography{refs2}

\end{document}